\newif\ifTwoColumn
\newif\ifTechReport

\TechReporttrue    

\ifTwoColumn
	\documentclass[twocolumn]{autart}

\else
	\ifTechReport
		\documentclass[eview,3p,number,11pt]{elsarticle_arXiv}
		\usepackage{amsmath,amsthm}
		\newenvironment{pf}{\begin{proof}}{\end{proof}}
		\usepackage{parskip}
		\usepackage{geometry}
	\else
		\documentclass[onecolumn]{autart}
	\fi
\fi

\usepackage{graphicx}
\usepackage[dvipsnames]{xcolor}
\graphicspath{{./figures/}}

\usepackage{amsmath,amssymb,amsfonts,enumitem,mathtools,bbm,mathrsfs}
\usepackage{dsfont}
\usepackage{epstopdf}
\usepackage{subfigure}
\usepackage[acronym]{glossaries}
\usepackage{balance}
\usepackage{url}
\usepackage[ruled]{algorithm}
\usepackage{booktabs}
\usepackage[official]{eurosym}
\usepackage{nicefrac}

\newtheorem{theorem}{Theorem}
\newtheorem{definition}{Definition}
\newtheorem{proposition}{Proposition}
\newtheorem{lemma}{Lemma}

\newtheorem{example}{Example}
\newtheorem{remark}{Remark}
\newtheorem{standing}{Standing Assumption}
\newtheorem{assumption}[standing]{Assumption}

\newcommand{\expected}{\mathbb{E}}
\newcommand{\R}{\mathbb{R}}

\newcommand{\N}{\mathbb{N}}

\newcommand{\mc}[1]{\mathcal{#1}}

\newcommand{\rank}{\mathrm{rank}}
\newcommand{\col}{\textrm{col}}

\newcommand{\bs}{\boldsymbol}
\newcommand{\bsone}{\boldsymbol{1}}

\newcommand{\blue}{\ignorespaces}

\newcommand\oprocendsymbol{\hbox{$\square$}}
\newcommand\oprocend{\relax\ifmmode\else\unskip\hfill\fi\oprocendsymbol}

\usepackage{xcolor,calc}

\makeglossaries
\newacronym{NEP}{NEP}{Nash equilibrium problem}
\newacronym{GNEP}{GNEP}{generalized Nash equilibrium problem}
\newacronym{iid}{i.i.d.\@}{independent and identically distributed}
\newacronym{wrt}{w.r.t.\@}{with respect to}
\newacronym{od}{OD}{origin-destination}
\newacronym{ncp}{NCP}{nonlinear complementarity problem}
\newacronym{LTI}{LTI}{linear time-invariant}
\newacronym{PEV}{PEV}{plug-in electric vehicle}
\newacronym{SoC}{SoC}{State of Charge}

\newglossaryentry{LP}
{
	name={LP},
	description={linear programming},
	first={\glsentrydesc{LP} (\glsentrytext{LP})},
	plural={LPs},
	descriptionplural={linear programs},
	firstplural={\glsentrydescplural{LP} (LPs)}
}

\newglossaryentry{QVI}
{
	name={VI},
	description={quasi variational inequality},
	first={\glsentrydesc{QVI} (\glsentrytext{QVI})},
	plural={VIs},
	descriptionplural={quasi variational inequalities},
	firstplural={\glsentrydescplural{QVI} (QVIs)}
}

\newglossaryentry{VI}
{
	name={VI},
	description={variational inequality},
	first={\glsentrydesc{VI} (\glsentrytext{VI})},
	plural={VIs},
	descriptionplural={variational inequalities},
	firstplural={\glsentrydescplural{VI} (VIs)}
}

\newglossaryentry{v-GNE}
{
	name={v-GNE},
	description={variational generalized Nash equilibrium},
	first={\glsentrydesc{v-GNE} (\glsentrytext{v-GNE})},
	plural={v-GNE},
	descriptionplural={variational generalized Nash equilibria},
	firstplural={\glsentrydescplural{v-GNE} (\glsentryplural{v-GNE})}
}

\newglossaryentry{GNE}
{
	name={GNE},
	description={generalized Nash equilibrium},
	first={\glsentrydesc{GNE} (\glsentrytext{GNE})},
	plural={GNE},
	descriptionplural={generalized Nash equilibria},
	firstplural={\glsentrydescplural{GNE} (\glsentryplural{GNE})}
}

\begin{document}

\begin{frontmatter}

\title{Probabilistic feasibility guarantees for solution sets to uncertain variational inequalities}

\author{Filippo Fabiani, Kostas Margellos and Paul J. Goulart}

\thanks{\emph{Email addresses:} \texttt{{filippo.fabiani, kostas.margellos, paul.goulart}@eng.ox.ac.uk}.
	This work was partially supported through the Government’s modern industrial strategy by Innovate UK, part of UK Research and Innovation, under Project LEO (Ref. 104781).}

\begin{abstract}
	We develop a data-driven approach to the computation of a-posteriori feasibility certificates to the solution sets of variational inequalities affected by uncertainty. Specifically, we focus on instances of variational inequalities with a deterministic mapping and an uncertain feasibility set, and represent uncertainty by means of scenarios. Building upon recent advances in the scenario approach literature, we quantify the robustness properties of the entire set of solutions of a variational inequality, with feasibility set constructed using the scenario approach, against a new unseen realization of the uncertainty. \blue{Our results extend existing results that typically impose an assumption that the solution set is a singleton and require certain non-degeneracy properties, and thereby offer probabilistic feasibility guarantees to any feasible solution.} We show that assessing the violation probability of an entire set of solutions, rather than of a singleton,  requires enumeration of the support constraints that ``shape'' this set. Additionally, we propose a general procedure to enumerate the support constraints that does not require a closed form description of the solution set, which is unlikely to be available. We show that robust game theory problems can be modelling via uncertain variational inequalities, and illustrate our theoretical results through extensive numerical simulations on a case study involving an electric vehicle charging coordination problem.
\end{abstract}

\end{frontmatter}


\section{Introduction}

As a general purpose tool embracing a rich class of decision-making problems, \glspl{VI} have been widely adopted in many scientific areas, from operation research and mathematical programming to optimization and game theory \cite{giannessi2006equilibrium,facchinei2007finite}.
Formally, a \gls{VI} is defined by means of a feasibility set $\mc{X} \subseteq \R^n$, and a mapping $F : \mc{X} \to \R^n$. We denote by VI$(\mc{X}, F)$ the problem of finding some vector $x^\star \in \mc{X}$ such that
\begin{equation}\label{eq:VI}
(y - x^\star)^\top F(x^\star) \geq 0, \, \text{ for all } y \in \mc{X}.
\end{equation}
Finding a point $x^\star$ satisfying \eqref{eq:VI} amounts to solving a generalized nonlinear complementarity problem, and therefore encompasses a broad variety of equilibrium problems that appear in multiple engineering domains \cite{ferris1997engineering}. Prominent examples include network and traffic problems \cite{dafermos1980traffic,dupuis1993dynamical,giannessi1995variational}, optimal control \cite{mignot1984optimal,friedman1986optimal,eckstein1998operator}, economics and demand-side management \cite{zhao1991general,daniele2003evolutionary,nagurney2005supply}.

However, most of the analytical and algorithmic results  in the literature on \glspl{VI} are restricted to deterministic problems.  Hence all of the aforementioned applications inherently neglect potential sources of uncertainty that strongly affect the data of the problem itself, i.e., the mapping $F$ and the feasible set $\mc{X}$. To move beyond the limited scope of deterministic methods, we focus on stochastic approaches to uncertain \glspl{VI}.   The existing literature in this area is broadly split into two main directions for incorporating uncertainty into the model in \eqref{eq:VI} \cite{shanbhag2013stochastic}: an expectation-based and a worst-case formulation.

Specifically, given a random variable $\delta \in \Delta$, we refer to the expected-value formulation of an uncertain \gls{VI}, originally described in \cite{king1993asymptotic}, as the problem of computing a deterministic vector $x^\star \in \mc{X}$ such that
\begin{equation}\label{eq:VI_expected}
	(y - x^\star)^\top \expected[F(x^\star, \delta)] \geq 0, \, \text{ for all } y \in \mc{X},
\end{equation}
where, in this case, $F : \mc{X} \times \Delta \to \R^n$. Unfortunately, there are only a few cases in which the formulation in \eqref{eq:VI_expected} is known to be computationally tractable, e.g., when $\delta$ takes values in a discrete set. In most stochastic regimes, the computation of the expected value requires multidimensional integration, which is a difficult task.

Moreover, for applications in which the decision arising from the solution to an uncertain \gls{VI} is required to be robust to parametric uncertainties, the formulation in \eqref{eq:VI_expected} may be inappropriate.  In those cases a worst-case formulation may be more suitable, allowing one to incorporate both uncertainty in the mapping $F(\cdot)$ and in the constraints $\mc{X}$. Specifically, the feasible set can be modelled as the (infinite) intersection of sets $\mc{X}_\delta$ generated by every possible realization of the random variable $\delta$, namely $\cap_{\delta \in \Delta} \mc{X}_\delta$. Thus, a vector $x^\star \in \cap_{\delta \in \Delta} \mc{X}_\delta$ is considered a solution to the uncertain \gls{VI} if
\begin{equation}\label{eq:VI_almost}
(y - x^\star)^\top F(x^\star, \delta) \geq 0, \, \text{ for all } y \in \left.\bigcap_{\delta \in \Delta} \mc{X}_\delta\right., \, \delta \in \Delta.
\end{equation}
However, such a worst-case formulation imposes two main challenges: i) the set $\Delta$ may be unknown and the only information available may come via data/scenarios for $\delta$; ii) even if $\Delta$ is known, it might be a continuous set (more generally, a set with infinite cardinality), thereby giving rise to an infinite set of constraints in \eqref{eq:VI_almost}.

To address these challenges, we \blue{adopt the data-driven approach proposed in \cite{campi2018general}} to quantify a-posteriori the feasibility of the entire set of solutions to the \gls{VI} against previously unseen realizations of the uncertainty. Specifically, in this paper we focus on a subclass of uncertain \glspl{VI} defined as in \eqref{eq:VI_almost} -- namely those instances characterized by a deterministic mapping $F(\cdot)$ and an uncertain feasible set $\mc{X}$. \blue{Using a set-oriented perspective, we recast our problem to the form of the abstract decision-making problems considered in \cite{campi2018general}. This enables us to inherit the probabilistic feasibility results established in \cite[Th.~1]{campi2018general}, and thereby characterize the robustness properties of the entire solution set to an uncertain~\gls{VI}.}

Remarkably, the family of \glspl{VI} we \blue{investigate} is quite large, encompassing applications across several domains:
\begin{itemize}[leftmargin=*]
\item \blue{ A wide class of \glspl{NEP} and \glspl{GNEP} can be characterized} as \glspl{VI} \cite{facchinei2007finite,facchinei2007generalized} \blue{in several applications of practical interest}. \blue{We show that the robust variant of \glspl{NEP}/\glspl{GNEP}} falls directly into the class of \glspl{VI} investigated in this paper \cite{aghassi2006robust,pantazis2020aposteriori,fabiani2020scenario}.

\item  As a static assignment problem, an optimal network flow in traffic networks can also be computed via solution to an associated \gls{VI} \cite[Th.~3.14]{patriksson1994traffic}. \blue{In this domain of application}, it is quite common to model the overall traffic demand as an uncertain variable \cite{gwinner2006random,siu2008doubly,chen2012vulnerability,daniele2015random,jadamba2018efficiency}, thus randomly constraining the traffic flow over admissible paths of the network.
\item Finite horizon control problems are typically formalized as constrained optimization problems, which can be modelled as \glspl{VI} \cite[\S 1.3.1]{facchinei2007finite}. Specifically, the optimal sequence of control inputs is required to minimize some predefined cost function, while being subject to operational and dynamic constraints, both of which may be affected by uncertainty \cite{rockafellar1986lagrangian,rockafellar1990generalized,zhang2017robust}.
\end{itemize}

	\subsection{Literature review and main contributions}
	To the best of our knowledge, this work is the first to address the problem of evaluating the robustness of the entire set of solutions to an uncertain \gls{VI} in a distribution-free fashion. Compared to existing results on data-driven approaches to assess the robustness of solutions to general \glspl{VI} (or to particular cases thereof), we consider a broad family of uncertain \glspl{VI} in \eqref{eq:VI_almost} rather than just \gls{VI} problems arising from the computation of \glspl{v-GNE}, a subset of \glspl{GNE} in \glspl{GNEP} \cite{fele2019probabilistic,fele2019probably,pantazis2020aposteriori,fabiani2020scenario}.

	Specifically, a \gls{NEP} is considered in \cite{fele2019probabilistic,fele2019probably} where the uncertain parameter is encoded as a common term in the agents’ cost function, while the constraint set of each player is deterministic. In this context, \cite{fele2019probabilistic,fele2019probably} provide an a-posteriori certificate on the probability that a (non-unique) variational solution to the \gls{NEP} remains unaltered upon a new realization of the uncertainty. The present work is complementary to the one in \cite{fele2019probabilistic,fele2019probably}, as they (indirectly) investigate \glspl{VI} with an uncertain mapping and a deterministic feasible set.

	Conversely, a \gls{NEP} with uncertain, yet affine, local constraints is considered in \cite{pantazis2020aposteriori}. By assuming deterministic cost functions, a contribution of \cite{pantazis2020aposteriori} is to provide robustness certificates for the constraint violation of any feasible point (thus including variational equilibria as special case) of the game considered. In contrast, we show in \S 3 that assessing the robustness of an equilibrium at a point inside the feasible set may lead to an over-conservative bound compared to the one derived in this paper, which is tailored for the entire set of equilibria.

	The approach proposed in \cite{fabiani2020scenario}, instead, paves the way to the set-oriented perspective investigated in this paper. Specifically, \cite{fabiani2020scenario} leverages the specific structure of the game in question, i.e., a \gls{GNEP} in aggregative form, to design probabilistic bounds on the feasibility of the entire set of variational equilibria.

	Finally, \cite{paccagnan2019scenario} has addressed robustness questions for uncertain \glspl{VI}, providing a-posteriori robustness certificates for the solution to uncertain (quasi-)\gls{VI} in \eqref{eq:VI_almost}.
	However, the perspective and proof line is substantially different from the one adopted in this paper. Specifically, in \cite{paccagnan2019scenario} it is postulated that the \gls{VI} admits a \emph{unique} solution, while certain non-degeneracy assumptions are imposed. Uniqueness restricts the class of \glspl{VI} that can be captured by problems of the type \eqref{eq:VI_almost}, while non-degeneracy is in general hard to verify even in optimization problems, and even moreso in \glspl{VI} and games \cite{campi2018wait,fele2019probably}.
	By considering a specific instance of the family of uncertain \glspl{VI} in \eqref{eq:VI_almost}, focusing on the entire set of solutions allow us to bypass both  assumptions.

	To conclude, we summarize our main contributions as follows:
	\begin{itemize}
		\item We consider a broad family of uncertain VIs in \eqref{eq:VI_almost} rather than just \gls{VI} problems arising in computing \gls{v-GNE}, thus complementing the results in \cite{pantazis2020aposteriori,fabiani2020scenario};
		\item In the spirit of \cite{fabiani2020scenario}, by relying on the data-driven tools given in \cite{campi2018general} we provide a-posteriori robustness certificates for the entire set of solutions to an uncertain \glspl{VI};
		Note that our set-oriented perspective is crucial for two reasons:
		\begin{enumerate}
			\item We are able to bypass the uniqueness and non-degeneracy assumptions postulated in \cite{paccagnan2019scenario};
			\item Compared to \cite{pantazis2020aposteriori}, we show that our bounds are, in general, less conservative;
			\item We offer guarantees for any feasible solution; hence we can support possibly suboptimal solutions returned by a generic algorithm.
		\end{enumerate}
		\item Our robustness certificates depend strongly  on the number of support subsamples characterizing the set of solutions to the uncertain \gls{VI}. We show that computing these support subsamples requires only an enumeration of the constraints that ``shape'' the solution set. An explicit representation of the unknown set of solutions is therefore not needed. In the case of affine constraints, we design a procedure to compute these samples that, in general, requires fewer iterations compared to the one in \cite{paccagnan2019scenario,campi2018general}.
	\end{itemize}
Finally, we show that problems in robust game theory falls within the class of uncertain \glspl{VI} we consider. Our theoretical results are supported through extensive numerical simulations on a \gls{GNEP} modelling the charging coordination of a fleet of \glspl{PEV}.

\subsection{Paper organization}
We formalize the data-driven problem addressed and state the main result of the paper, i.e., Theorem~\ref{th:VI}, in \S 2. In \S 3 we  discuss how the set-oriented problem we consider can be recast in the framework proposed in \cite{campi2018general}, thereby paving the way for a formal proof of Theorem~\ref{th:VI}. We then describe a systematic procedure to compute the number of support subsamples in the case of affine constraints in \S 4, also discussing the computational aspects associated with the proposed approach. Finally, in \S 5 we demonstrate our theoretical results through a numerical simulations on a \gls{GNEP}.

\subsection{Notation}
	\textit{Basic notation}: $\N$, $\R$, and $\R_{\geq 0}$ denote the set of natural, real, and nonnegative real numbers, respectively, with $\N_0 \coloneqq \N \cup \{0\}$. Given some $x \in \R^n$, $\|x\|$ is the Euclidean norm. Denote vectors of appropriate dimensions with elements all equal to $1$ ($0$) as $\bs{1}$ ($\bs{0}$). Given a matrix $A \in \R^{m \times n}$, its $(i,j)$ entry is denoted by $a_{i,j}$, $A^\top$ denotes its transpose, while for $A \in \R^{n \times n}$, $A \succ 0$ ($\succcurlyeq 0$) implies that $A$ is symmetric and positive (semi)-definite. For $A \succ 0$, $\| x \|_A := \sqrt{ x^\top A x }$. $\mc{C}^1$ is the class of continuously differentiable functions.
	For a given set $\mc{S} \subseteq \R^n$, $|\mc{S}|$ represents its cardinality, and $\textrm{int}(\mc{S})$, $\textrm{relint}(\mc{S})$ and $\textrm{bdry}(\mc{S})$ denote its topological interior, relative interior and boundary, respectively.  The set $\textrm{aff}(\mc{S})$ denotes its affine hull, i.e., the smallest affine set containing $\mc{S}$. The operator $\otimes$ denotes the Kronecker product, while $\col(\cdot)$ stacks its arguments in column vectors or matrices of compatible dimensions.
	For vectors $v_1,\dots,v_N\in\mathbb{R}^n$ and $\mc I=\{1,\dots,N \}$, we denote $\boldsymbol{v} \coloneqq (v_1 ^\top,\dots ,v_N^\top )^\top = \mathrm{col}((v_i)_{i\in\mc I})$ and $ \bs{v}_{-i} \coloneqq \col(( v_j )_{j\in\mc I\setminus \{i\}})$. With a slight abuse of notation, we sometimes use $\bs{v} = (v_i,\bs{v}_{-i})$.

	\textit{Operator-theoretic definitions}: Given a function $\phi : \R^n \to \R$, $\textrm{dom}(\phi) \coloneqq \{x \in \R^n \mid \phi(x) < \infty\}$ is the domain of $\phi$; $\partial \phi : \textrm{dom}(\phi) \rightrightarrows \R^n$ denotes the subdifferential set-valued mapping of $\phi$, defined as $\partial \phi(x) \coloneqq \{d \in \R^n \mid \phi(z) \geq \phi(x) + d^\top(z - x), \forall z \in \textrm{dom}(\phi)\}$, for all $x \in \textrm{dom}(\phi)$. For a given set $\mc{S} \subseteq \R^n$, the mapping $T: \mc{X} \rightarrow \R^n$ is pseudomonotone on $\mc{S}$ if for all $x$, $y \in \mc{S}$, $(x - y)^\top T(y) \geq 0 \implies (x - y)^\top T(x) \geq 0$; (strictly) monotone if $(T(x) - T(y))^\top(x - y) \, (>) \geq  0$ for all $x, y \in \mc{S}$ (and $x \neq y$); strongly monotone if there exists a constant $c > 0$ such that $(T(x) - T(y))^\top(x - y) \geq c \|x - y\|^2$ for all $x, y \in \mc{S}$.
	If $\mc{S}$ is nonempty and convex, the normal cone of $\mc{S}$ evaluated at  $x$ is the set-valued mapping $\mc{N}_{\mc{S}} : \R^n \rightrightarrows \R^n$, defined as $\mc{N}_{\mc{S}}(x) \coloneqq \{	d \in \R^n \mid d^\top (y - x) \leq 0, \; \forall y \in \mc{S}	\}$ if $x \in \mc{S}$, $\mc{N}_{\mc{S}}(x) \coloneqq \emptyset$ otherwise.

\section{Problem statement and main result}
We start by recalling some key concepts about \glsfirstplural{VI} \cite{facchinei2007finite}, and then describe the data-driven problem we consider throughout the paper.  We also state the main result of the paper\blue{, i.e., Theorem~\ref{th:VI}, in this section, but will defer the proof of this	 result to \S \ref{sec:scenario_VI}.} Unless otherwise specified, we assume measurability of all the quantities introduced hereafter.

\subsection{Background on \glspl{VI}}
Let us consider the deterministic \gls{VI} formally introduced in \eqref{eq:VI}, and let $\Omega \subseteq \mc{X}$ be the set of solutions to VI$(\mc{X}, F)$. The relation in \eqref{eq:VI} has a strong geometric interpretation that relies on the definition of the normal cone \blue{\cite[Ch.~1.1]{facchinei2007finite}}. If $\mc{X}$ is nonempty and convex, a vector $x^\star \in \mc{X}$ solves VI$(\mc{X}, F)$ if and only if $-F(x^\star) \in \mc{N}_\mc{X}(x^\star)$. For example, in the specific case of an optimization problem $\textrm{min}_{x \in \mc{X}} \, f(x)$, we have $F = \nabla f$ and the inclusion $-\nabla f (x^\star) \in \mc{N}_\mc{X}(x^\star)$ corresponds to satisfaction of the KKT conditions at some $x^\star$. Thus, in view of the definition of the normal cone, any point belonging to $\textrm{int}(\mc{X})$ solves VI$(\mc{X}, F)$ if and only if $F(x^\star) = 0$.

The structural properties of both the feasible set $\mc{X}$ and the mapping $F(\cdot)$ allow one to establish the existence and uniqueness of the solution to VI$(\mc{X}, F)$, as well as to provide the minimal conditions that enable one to design suitable solution algorithms with convergence guarantees. By combining \cite[Cor.~2.2.5, Th.~2.3.5]{facchinei2007finite}, we can characterize the solution set to VI$(\mc{X}, F)$ as follows:
\begin{lemma}\label{lemma:VI_sol_set}
	Let $\mc{X}$ be a compact and convex set, and let $F(\cdot)$ be a continuous mapping. Then, the following statements hold true:
	\begin{itemize}
		\item[$\mathrm{(i)}$] $\Omega$ is a nonempty and compact set;
		\item[$\mathrm{(ii)}$] If $F(\cdot)$ is also pseudomonotone, then $\Omega$ is also a convex set.
	\end{itemize}
\end{lemma}
Note that assuming strong monotonicity of $F(\cdot)$ would guarantee the existence of a unique solution to VI$(\mc{X}, F)$ \cite[Th.~2.3.3]{facchinei2007finite}. Requiring the mapping $F(\cdot)$ to be only pseudomonotone is clearly weaker than assuming (strong) monotonicity. However, pseudomonotonicity is not always a trivial condition to verify, while monotonicity is naturally satisfied in several practical applications that involve, e.g., the subdifferential $\partial f$ of a proper, closed, convex function $f : \R^n \to \R$ (or its conjugate, see \cite[\S 20]{bauschke2011convex}). \blue{Both conditions, however, individually represent} one of the weakest assumptions that guarantee convergence of many efficient solution algorithms, see, e.g., \blue{\cite{el2001pseudomonotone}, \cite{he2002new}, \cite{malitsky2015projected}, \cite{nguyen2018extragradient}, \cite{vuong2018weak}, or the dedicated sections in} \cite[\S 7, \S 12]{facchinei2007finite}, \cite[Ch.~12]{palomar2010convex} \blue{and references therein}.

\subsection{Uncertain \glspl{VI} and scenario-based formulation}\label{subsec:unc_VI}
We aim to provide out-of-sample feasibility certificates for the entire set of solutions to a given uncertain \gls{VI} by exploiting some observed realizations, i.e., \emph{scenarios}, of the uncertain parameter $\delta$. Formally, let us consider a probability space $(\Delta,\mc{D},\mathbb{P})$, where $\Delta \subseteq \R^\ell$ represents the set of values that $\delta$ can take, $\mc{D}$ is a $\sigma$-algebra and $\mathbb{P}$ is a (possibly unknown) probability measure over $\mc{D}$. Given a mapping $F : \mc{X} \rightarrow \R^n$ and a deterministic feasible set $\mc{X} \subseteq \R^n$, let $\mc{X}_\delta \subseteq \R^n$ be an additional set of constraints associated with the uncertain parameter $\delta$. We define the worst-case \gls{VI} problem, denoted as VI$(\mc{X} \cap \mc{X}_\delta, F)$, as the problem of finding some $x^\star \in \mc{X} \cap \mc{X}_\delta$ that satisfies
\begin{equation}\label{eq:unc_VI}
(y - x^\star)^\top F(x^\star) \geq 0, \text{ for all } y \in \mc{X} \cap \mc{X}_\delta, \; \delta \in \Delta.
\end{equation}
However, given the possibly infinite cardinality of $\Delta$ in \eqref{eq:unc_VI} and motivated by the increasing availability of data, we investigate a data-driven approach. Specifically, let us consider $\delta_K \coloneqq \{\delta^{(i)}\}_{i \in \mc{K}} = \{\delta^{(1)}, \ldots, \delta^{(K)}\} \in \Delta^K$, $\mc{K} \coloneqq \{1,2,\ldots,K\}$, hereafter also called a $K$-multisample, as a finite collection of $K \in \N$ \gls{iid} observed realizations of $\delta$. Here, every $K$-multisample is defined over the probability space $(\Delta^K, \mc{D}^K, \mathbb{P}^K)$, resulting from the $K$-fold Cartesian product of the original probability space $(\Delta,\mc{D},\mathbb{P})$.
Let $\mc{X}_{\delta^{(i)}}$ be a constraint set associated with the $i$-th sample, which constrains the decisions that are admissible for the situation represented by $\delta^{(i)}$.
The scenario-based \gls{VI} problem VI$(\mc{X}_{\delta_K}, F)$, \blue{with $\mc{X}_{\delta_K} \coloneqq \cap_{i \in \mc{K}} \mc{X}_{\delta^{(i)}} \cap \mc{X}$}, is then the problem of finding an $x^\star \in \mc{X}_{\delta_K}$ such that
\begin{equation}\label{eq:scenario_based_VI}
(y - x^\star)^\top F(x^\star) \geq 0, \text{ for all } y \in \mc{X}_{\delta_K}.
\end{equation}
Let us define the set of solutions to \eqref{eq:scenario_based_VI} as
\begin{equation}\label{eq:sol_set}
\Omega_{\delta_K} \coloneqq \{x \in \mc{X}_{\delta_K} \mid (y - x)^\top F(x) \geq 0, \, \forall y \in \mc{X}_{\delta_K}\}.
\end{equation}
Note that given the dependency on the set of $K$ realizations $\delta_K$, the set $\Omega_{\delta_K}$ is itself a random variable. For the case $K = 0$ our problem reduces to a deterministic \gls{VI} problem, i.e., VI$(\mc{X},F)$, where no uncertainty is present. Let $\Omega_{\delta_0}$ be the solution set for this case.
In light of the results of Lemma~\ref{lemma:VI_sol_set}, we will make the following assumptions throughout the remainder of the paper:
\begin{standing}\label{standing:sets}
	For any $K \in \N_0$, the set $\mc{X}_{\delta_K}$ is a nonempty, compact and convex set for all $\delta_K \in \Delta^K$.
\end{standing}
\begin{standing}\label{standing:mapping}
	The mapping $F : \mc{X} \to \R^n$ is continuous and pseudomonotone.
\end{standing}

These assumptions ensure that our scenario-based \gls{VI} \eqref{eq:scenario_based_VI} has a non-empty solution set:

\begin{lemma}\label{lemma:solution_compact}
	For all $K \in \N_0$, $\Omega_{\delta_K}$ is a nonempty, compact and convex set.
\end{lemma}
\begin{pf}
	It follows immediately from Lemma~\ref{lemma:VI_sol_set}(ii), as $F(\cdot)$ is continuous and pseudomonotone and $\mc{X}_{\delta_K}$ is a finite intersection (due to Standing Assumption~\ref{standing:sets}) of nonempty, compact and convex sets.
\end{pf}
In the spirit of \cite{campi2018general}, we introduce $\Theta_K : \Delta^K \rightrightarrows \mc{X}$ as the mapping that, given a set of realizations $\delta_K$, returns the solution set to VI$(\mc{X}_{\delta_K},F)$, namely 
\begin{equation}\label{eq:theta_mapping}
	\Theta_K(\delta^{(1)}, \ldots, \delta^{(K)}) = \Theta_K(\delta_K) \coloneqq \Omega_{\delta_K} .
\end{equation}
When $K = 0$, we assume that $\Theta_0$ returns the solution set to VI$(\mc{X}, F)$, $\Omega_{\delta_0}$.


\subsection{Robustness certificates for solution sets to \glspl{VI}}
Given any $K$-multisample $\delta_K$, we are interested in evaluating the robustness of the entire set of solutions $\Omega_{\delta_K}$ in \eqref{eq:sol_set} to a previously unseen realization of the uncertain parameter $\delta$. To this end, let $\Omega_\delta$ be the set of solutions induced by a certain realization $\delta \in \Delta$. We now introduce the following definition of violation probability of a \blue{generic} set \blue{of solutions, $\Omega$}.
\begin{definition}\textup{(Violation Probability of a Set)}\label{def:violation_set}
	The \emph{violation probability} associated with a set \blue{of solutions} $\Omega$ is defined as
	\begin{equation}\label{eq:violation_set}
	V(\Omega) \coloneqq \mathbb{P}\{\delta \in \Delta \mid \Omega \not\subseteq \Omega_\delta \}.
	\end{equation}
\end{definition}
Informally speaking, the condition $\Omega \not\subseteq \Omega_\delta$ implies that, once $\delta$ is drawn, at least one element in $\Omega$ ceases to be a solution.
Note that the set $\Omega_{\delta_K}$ is itself a random variable, and hence so is the violation probability $V(\Omega_{\delta_K})$. We therefore wish to characterise our confidence that $V(\Omega_{\delta_K})$ is below some violation level. Before stating the main result of this section, we recall the following definition that will be crucial for the remainder of the paper:
\begin{definition}\textup{(Support Subsample) \cite[Def.~2]{campi2018general}}\label{def:support_sub}
	Given any $\delta_K \in \Delta^K$, a \emph{support subsample} $S \subseteq \delta_K$ is a $p$-tuple of unique elements of $\delta_K$, i.e., $S \coloneqq \{\delta^{(i_1)}, \ldots, \delta^{(i_p)}\}$, $i_1 < \ldots < i_p$, that gives the same solution as the original sample, i.e.,
	$$
	\Theta_p(\delta^{(i_1)}, \ldots, \delta^{(i_p)}) = \Theta_K(\delta^{(1)}, \ldots, \delta^{(K)}).
	$$
\end{definition}
Here, let $\Upsilon_K : \delta_K \rightrightarrows \mc{K}$ be any algorithm returning a $p$-tuple $\{i_1, \ldots, i_p\}$, $i_1 < \ldots < i_p$,  such that $\{\delta^{(i_1)}, \ldots, \delta^{(i_p)}\}$ is a support subsample for $\delta_K$, and let $s_K \coloneqq |\Upsilon_K(\delta_K)|$.
Note that $s_K$ is itself a random variable since it depends on $\delta_K$. We discuss the construction of such an algorithm in \S 4. \blue{Our main result characterizes} the violation probability of $\Omega_{\delta_K}$, i.e., the solution set to the scenario-based \gls{VI} in \eqref{eq:scenario_based_VI}, as follows:

\begin{theorem}\label{th:VI}
	Fix $\beta \in (0,1)$, and let $\varepsilon : \mc{K} \cup \{0\} \to [0, 1]$ be a function such that
	\begin{equation}\label{eq:epsilon}
	\left\{\begin{aligned}
	& \varepsilon(K) = 1,\\
	& \sum_{h = 0}^{K - 1} \left( \begin{array}{c}
	K\\
	h
	\end{array} \right) (1 - \varepsilon(h))^{K - h} = \beta.
	\end{aligned}
	\right.
	\end{equation}
	Then, for any mappings $\Theta_K$, $\Upsilon_K$ and distribution $\mathbb{P}$, it holds that
	\begin{equation}\label{eq:prob_feas_boud}
	\mathbb{P}^K \{\delta_K \in \Delta^K \mid V(\Omega_{\delta_K}) > \varepsilon(s_K) \} \leq \beta.
	\end{equation}
\end{theorem}

Note that the bound in \eqref{eq:prob_feas_boud} is an a-posteriori statement since $s_K$ depends on the multisample extracted.  In words, Theorem~\ref{th:VI} implies that the probability that $\Omega_{\delta_K \cup \{\delta\}}$ differs from $\Omega_{\delta_K}$ (as $\Omega_{\delta_K} \subseteq \Omega_{\delta_K \cup \{\delta\}}$ necessarily implies that $\Omega_{\delta_K} = \Omega_{\delta_K \cup \{\delta\}}$ -- see also Lemma~\ref{lemma:if_affine}) is at most equal to $\varepsilon(s_K)$, with confidence at least $1-\beta$, for an arbitrarily small $\beta \in (0,1)$. We give the proof of Theorem~\ref{th:VI} in the next section, after first stating and proving some ancillary results.

\section{The scenario approach to uncertain \glspl{VI}}\label{sec:scenario_VI}
In this section, we first recall some key notions of the scenario approach theory, and then we show how they can be extended to the context of solution sets to \glspl{VI}. We finally conclude by proving and discussing Theorem~\ref{th:VI}.
\subsection{Scenario approach for decision-making problems}
The scenario approach was initially conceived to provide a-priori out-of-sample feasibility guarantees associated with the solution to an uncertain convex optimization problem \cite{calafiore2006scenario,calafiore2010random,campi2018introduction}.  It has recently been extended to abstract decision-making problems through an a-posteriori assessment of the feasibility risk  \cite{campi2018general,campi2018wait}.

With a slight abuse of notation, we assume here that $\Theta_K : \Delta^K \to \mc{X}$ represents a function leading to a single scenario decision $\theta^\star_{\delta_K}$ for some generic abstract decision-making problem, rather than as a set of solutions specific to a \gls{VI} as in \eqref{eq:sol_set}--\eqref{eq:theta_mapping}.
In accordance with \cite{campi2018general}, $\theta_{\delta_K}^\star$ is assumed to be unique, otherwise any convex tie-break rule may be employed \cite{campi2008exact}. Then, we recall the following assumption that is crucial to prove \cite[Th.~1]{campi2018general}.
\begin{assumption}\textup{\cite[Ass.~1]{campi2018general}}\label{ass_campi_garatti_ramponi}
	For all $K \in \N$ and for all $\delta_K \in \Delta^K$, it holds that $\Theta_K(\delta_K) \in \mc{X}_{\delta^{(i)}}$, for all $i \in \mc{K}$.
\end{assumption}
%
Assumption~\ref{ass_campi_garatti_ramponi} implies that the decision taken while observing $K$ realizations of the uncertainty $\delta$ is consistent \gls{wrt} all the extracted scenarios. The goal in \cite{campi2018general} was then to assess the violation probability of the scenario decision $\theta^\star_{\delta_K}$, as formalized next.
\begin{definition}\textup{(Violation Probability of a Singleton)}\label{def:viol_prob}
	The violation probability of a decision $\theta \in \mc{X}$ is given by
	$$
		V(\theta) \coloneqq \mathbb{P}\{\delta \in \Delta \mid \theta \notin \mc{X}_\delta\}.
	$$
\end{definition}

%
%
%
%
%
%
%
%
%
%

Notice again the slight abuse of notation, where we use $V$ to denote both the violation of a singleton $\theta$ in Definition~\ref{def:viol_prob} and of a set in Definition~\ref{def:violation_set}\blue{, while $\Theta_K$ in this subsection returns an element of $\mc{X}$ (the solution $\theta^\star_{\delta_K}$) rather than a set as in \S \ref{subsec:unc_VI}}. The results in \cite{campi2018general} hold for generic decisions, as long as Assumption~\ref{ass_campi_garatti_ramponi} is satisfied. In the next subsection, we show how we can employ those results and adapt the sequence of inclusions in Assumption~\ref{ass_campi_garatti_ramponi} when our decision is a set.
With the set-oriented perspective introduced in \S 2, for the uncertain \gls{VI} in \eqref{eq:unc_VI} we let the admissible decision for the situation represented by $\delta$ coincide with the solution set $\Omega_\delta$, which is clearly a subset of the feasible set shaped by the uncertain parameter, i.e., the set $\mc{X} \cap \mc{X}_\delta$. This clarifies the analogy between Definition~\ref{def:violation_set} and \ref{def:viol_prob}.
For completeness, we restate \blue{as a lemma} the crucial result provided in \cite{campi2018general} to bound the violation probability of $\theta^\star_{\delta_K}$.
\begin{lemma}\textup{\cite[Th.~1]{campi2018general}}\label{lemma:campigaratti}
	Let Assumption~\ref{ass_campi_garatti_ramponi} hold true and fix $\beta \in (0,1)$. Let $\varepsilon : \mc{K} \cup \{0\} \to [0, 1]$ be the function defined in \eqref{eq:epsilon}.
	Then, for the $\Theta_K$ defined in this subsection, and for any $\Upsilon_K$ as defined below Definition~\ref{def:support_sub}, we have that
	$
		\mathbb{P}^K \{\delta_K \in \Delta^K \mid V(\theta^\star_{\delta_K}) > \varepsilon(s_K) \} \leq \beta.
	$
\end{lemma}

\blue{Notice that, in this case, $s_K$ would be the number of samples such that, by feeding $\Theta_K$ only with those samples, would return the same optimal solution $\theta^\star_{\delta_K}$ that would have been obtained if all samples were employed.}
We will use Lemma~\ref{lemma:campigaratti} to prove our main result in Theorem~\ref{th:VI}, which characterizes the entire set of solutions to VI$(\mc{X}_{\delta_K},F)$, but require some preliminary results first.

\subsection{Scenario-based \gls{VI} solution sets as nested sets of decisions}

%
%

In view of the analogy between Definition~\ref{def:violation_set} and \ref{def:viol_prob}, we aim to follow the approach of \cite{campi2018general} by focusing on a set of decisions, extending the conditions in \S 3.1, and in particular the sequence of inclusions in Assumption~\ref{ass_campi_garatti_ramponi}, to the solution set for the uncertain \gls{VI} in \eqref{eq:unc_VI}.
Thus, returning to the more general case where $\Theta_K$ is a set-valued mapping as defined in \eqref{eq:theta_mapping}, i.e.\ $\Theta_K : \Delta^K \rightrightarrows \mc{X}$, since we focus on the entire set of solutions, we remark that, for any $K$-multisample $\delta_K \in \Delta^K$, the uniqueness of the solution returned by $\Theta_K$ holds by definition.  Then, in the spirit of \cite[Def.~2]{margellos2015connection}, we envision that the set-oriented counterpart of the sequence of inclusions in Assumption~\ref{ass_campi_garatti_ramponi} shall be naturally translated into a consistency property of $\Omega_{\delta_K}$, as defined next.
\begin{definition}\textup{(Consistency of Solution Sets)}\label{def:cons_set}
	Given some $K \in \N$ and $\delta_K \in \Delta^K$, the solution set to $\mathrm{VI}(\mc{X}_{\delta_K}, F)$ is \emph{consistent} with the collected scenarios if $\Theta_K(\delta_K) = \Omega_{\delta_K} \subseteq \mc{X}_{\delta^{(i)}}$, for all $i \in \mc{K}$.
\end{definition}


\blue{In analogy with Assumption~\ref{ass_campi_garatti_ramponi}, Definition~\ref{def:cons_set} establishes that the set of solutions to $\mathrm{VI}(\mc{X}_{\delta_K}, F)$, $\Omega_{\delta_K}$, which is based on $K$ scenarios, should be feasible for each of the sets $\mc{X}_{\delta^{(i)}}$, $i \in \mc{K}$, corresponding to each of the $K$ realizations of the uncertain parameter.}
Thus, a first step towards applying the bound in \blue{Lemma~\ref{lemma:campigaratti}} is to show that the mapping $\Theta_K(\cdot)$ introduced in \eqref{eq:theta_mapping} is consistent with the realizations observed in the scenario-based \gls{VI} in \eqref{eq:scenario_based_VI}. This fact, however, follows by definition. For any $K \in \N$ and associated $K$-multisample $\delta_{K} \in \Delta^K$, indeed, it holds that $\Theta_K(\delta_K) \coloneqq \Omega_{\delta_K} \subseteq \cap_{i \in \mc{K}} \mc{X}_{\delta^{(i)}} \cap \mc{X}$, which on the other hand implies that $\Theta_K(\delta_K) \subseteq \mc{X}_{\delta^{(i)}}$, for all $i \in \mc{K}$, thus directly falling within Definition~\ref{def:cons_set}. We will make use of these considerations to rely on the bound in \blue{Lemma~\ref{lemma:campigaratti}} in the proof of Theorem~\ref{th:VI}, along with of the assumption on the solution set $\Omega_{\delta_K}$ introduced next.
\begin{standing}\label{standing:affinehull}
	For all $K \in \N$ and $\delta_K \in \Delta^K$, $\mathrm{aff}(\Omega_{\delta_K}) = \mathrm{aff}(\Omega_{\delta_0})$.
\end{standing}
Specifically, if the uncertain \gls{VI} in \eqref{eq:unc_VI} is defined in $\R^n$ and $\Omega_{\delta_K}$ is a convex, $m$-dimensional set, then Standing Assumption~\ref{standing:affinehull} allows for $m < n$. In this sense, assuming $\mathrm{aff}(\Omega_{\delta_K}) = \mathrm{aff}(\Omega_{\delta_0})$ for any $\delta_K \in \Delta^K$, $K \in \N$, is weaker than, e.g., assuming the nonemptiness of $\textrm{int}(\Omega_{\delta_K})$ for every possible realization of $\delta_K$. Standing Assumption~\ref{standing:affinehull} rules out the scenario that a given realization of the uncertainty $\delta_K$ reduces the solution set $\Omega_{\delta_K}$ to one of lower dimension with a different affine hull compared to the one of the deterministic \gls{VI}. To clarify its role, we introduce and discuss the following illustrative example.

\begin{figure}[t]
	\centering
	\ifTwoColumn
		\includegraphics[trim=1cm 0 .5cm 0,width=.9\columnwidth]{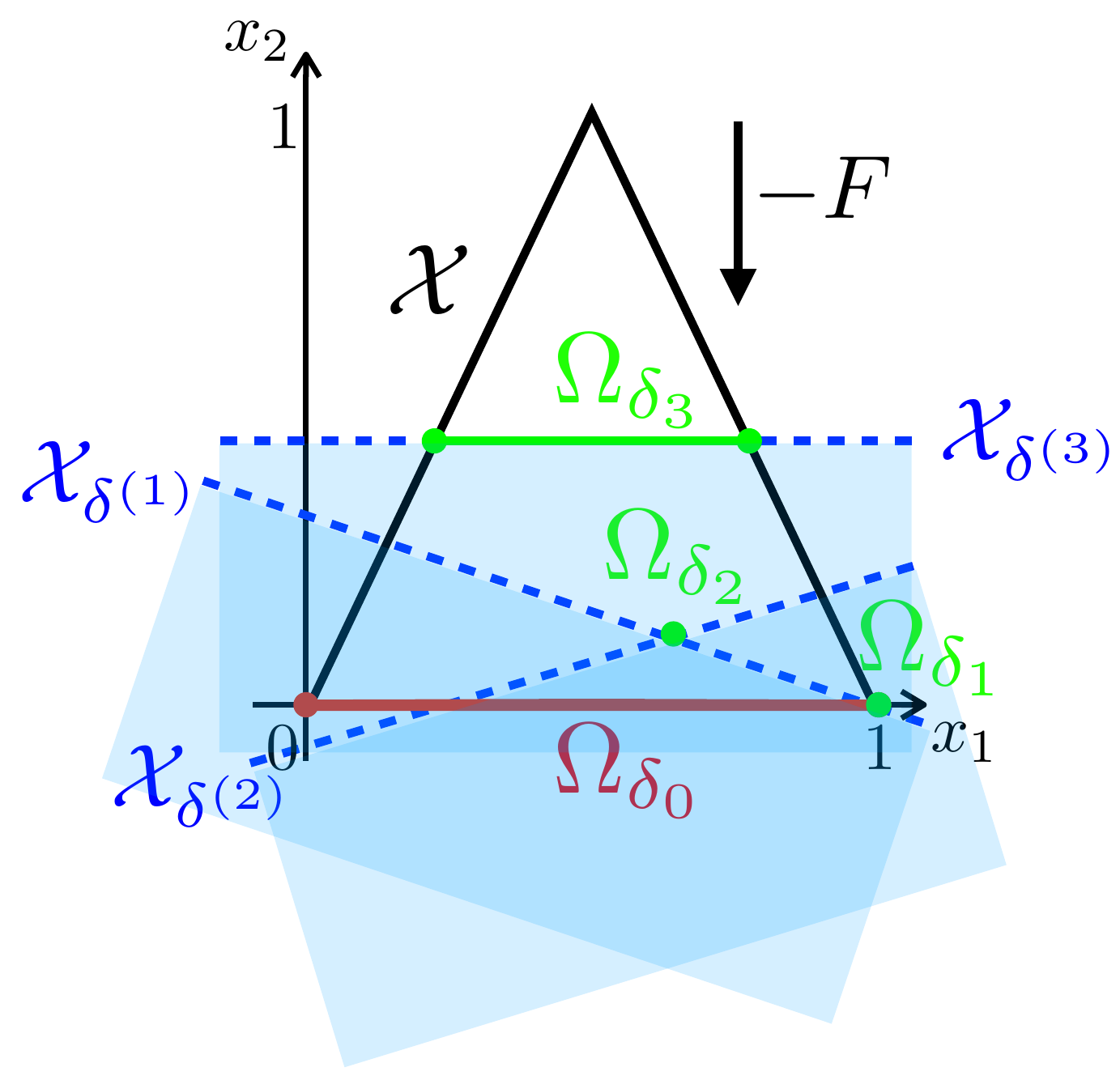}
	\else
	 	\includegraphics[trim=1cm 0 .5cm 0,width=0.45\columnwidth]{example_VI}
	\fi
	\caption{Compared to $\Omega_{\delta_0}$ (red line), every realization of $\delta$ (dashed blue lines, while the shaded cyan area denotes a region excluded by any $\mc{X}_{\delta^{(i)}}$, $i = 1, 2, 3$) results in a solution set $\Omega_{\delta_i}$, $i = 1, 2, 3$, that belongs to a different affine hull and/or on a space of lower dimension (green dots or line). Standing Assumption~\ref{standing:affinehull} allows us to rule out such cases.}
	\label{fig:example_VI}
\end{figure}
\begin{example}\label{ex:affine}
	Let us consider a two-dimensional case as shown in Fig.~\ref{fig:example_VI}, where $F = \mathrm{col}(0,-1)$, is monotone and $\mc{X}$ has a triangular shape. Here, $\Omega_{\delta_0} = \{ x \in \R^2 \mid x_1 \in [0,1], \, x_2 = 0 \}$, and its affine hull corresponds to the entire $x_1$-axis. After observing the first realization of $\delta$, i.e., $\delta^{(1)}$, which introduces the set $\mc{X}_{\delta^{(1)}} = \{ x \in \R^2 \mid -[\nicefrac{1}{3} \; 1]^\top \, x \leq -\nicefrac{1}{3} \}$, the solution set reduces to a singleton $\Omega_{\delta_1} = \{ x \in \R^2 \mid x = \mathrm{col}(1,0)\}$. Here, $\Omega_{1}$ has a smaller dimension compared to $\Omega_{\delta_0}$, despite its affine hull, i.e., the singleton itself, being a subset of the $x_1$-axis. Then, drawing a new sample $\delta^{(2)}$, which introduces the set $\mc{X}_{\delta^{(2)}} = \{ x \in \R^2 \mid [\nicefrac{1}{3} \; -1]^\top \, x \leq \nicefrac{1}{15}\}$, we have $\Omega_{\delta_2} = \{ x \in \R^2 \mid  x = \mathrm{col}(\nicefrac{3}{5},\nicefrac{2}{15}) \}$, which has the same dimension as $\Omega_{\delta_1}$ but its affine hull is not a subset of $\mathrm{aff}(\Omega_{\delta_0})$. Finally, the third sample, $\delta^{(3)}$, introduces the set $\mc{X}_{\delta^{(3)}} = \{ x \in \R^2 \mid [0 \; -1]^\top \, x \leq -\nicefrac{1}{2}\}$, and hence we have $\Omega_{\delta_3} = \{ x \in \R^2 \mid x_1 \in [\nicefrac{1}{4}, \nicefrac{3}{4}], \, x_2 = \nicefrac{1}{2} \}$. Here, $\Omega_{\delta_3}$ has the same dimension of $\Omega_{\delta_0}$ but a different affine hull, i.e., the $x_1$-axis translated to $x_2 = \nicefrac{1}{2}$. Standing Assumption~\ref{standing:affinehull} is meant to rule out all these possible scenarios, allowing only for samples that ``shape'' $\mathrm{aff}(\Omega_{\delta_0})$ without altering its dimension.
\end{example}
As we investigate uncertain \glspl{VI} of the form \eqref{eq:unc_VI} where the uncertainty $\delta$ affects the feasible set only, Example~\ref{ex:affine} provides insight on translating Standing Assumption~\ref{standing:affinehull} to a condition on the probability space $\Delta$. In fact, it represents situations that can generally happen with non-zero probability. Specifically, let $\Delta$ in Example~\ref{ex:affine} be a subset of $\R^2$, namely the uncertainty $\delta = \col(a, b)$ has two components, and let $a \in \R$ parametrize the slope and $b \in \R$ the offset of the halfspaces introduced by every scenario, i.e., $\mc{X}_{\delta} = \{ x \in \R^2 \mid [a \; 1] \, x \leq b \}$, for every $\delta \in \Delta$. Then, for any distribution that admits a density, we can find non-zero intervals for $a$ and $b$ such that \blue{the \gls{iid} scenarios $\delta$ can be extracted from a restricted subset of $\Delta$, determined by the values of $a$ and $b$ themselves,} in order to meet Standing Assumption~\ref{standing:affinehull}, and hence ruling out the \blue{pathological cases} shown in Example~\ref{ex:affine}.
\blue{In the case the samples are extracted from a restricted subset of $\Delta$, note that the guarantees would hold for the probability measure that is induced by this restriction, and not for the original uncertainty measure. Alternatively, if Standing Assumption~\ref{standing:affinehull} is not satisfied for all multisamples, then we can still claim that with confidence at most $\beta$, if Standing Assumption~\ref{standing:affinehull} is satisfied, then the probability of violation is greater than $\varepsilon(s_K)$. To achieve this, in the statement of Theorem~\ref{th:VI}, rather than $\Theta_K$ we can restrict the space of multisamples to the ones for which Standing Assumption~\ref{standing:affinehull} is satisfied. This is analogous to the way infeasible problem instances are accounted for in \cite{calafiore2010random,calafiore2006scenario}.}
Moreover, note that by adopting restrictions on $\Delta$ as described above, Standing Assumption~\ref{standing:affinehull} allows us to address the strongly monotone case, where VI$(\mc{X}\cap \mc{X}_\delta,F)$ has a unique solution, for all $\delta \in \Delta$.
\begin{remark}
	In view of \cite[Th.~2.4.15]{facchinei2007finite}, Standing Assumption~\ref{standing:affinehull} can be replaced by a simpler one, which is easier to verify, in all problems that involve a (monotone) affine $\mathrm{\gls{VI}}$ with polyhedral feasible set -- see, e.g., Assumption~2 in \cite{fabiani2020scenario}.
\end{remark}

\begin{figure}
	\centering
	\ifTwoColumn
		\includegraphics[width=.7\columnwidth,trim=1.5cm 1.5cm 1.5cm 1.5cm]{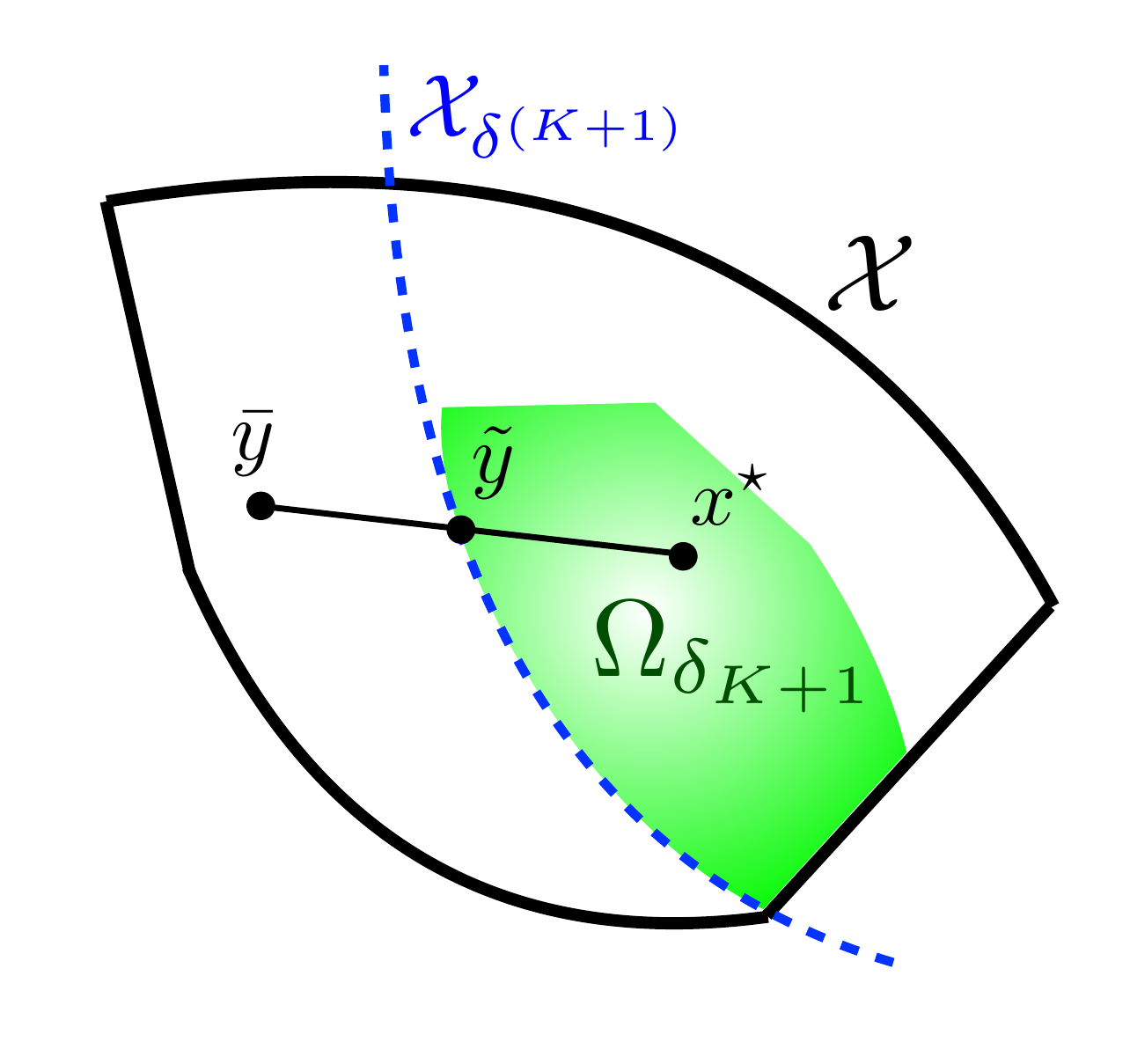}
	\else
	  \includegraphics[width=.35\columnwidth,trim=1.5cm 1.5cm 1.5cm 1.5cm]{proof_scheme.pdf}
	\fi
	\caption{Schematic two-dimensional construction of the proof of Lemma~\ref{lemma:if_affine}, part \textrm{(ii)}. Due to the convexity, there always exists some $\tilde{y} \in \mc{X}_{\delta_{K+1}}$, but $\tilde{y} \notin \textrm{int}(\Omega_{\delta_{K+1}})$, that allows to construct a contradiction. In this case, $\tilde{y} \in \textrm{bdry}(\Omega_{\delta_{K+1}})$.}\label{fig:proof_scheme}
\end{figure}

Given some $K \in \N$, let $\Omega_{\delta_{K+1}} \coloneqq \Omega_{\delta_K \cup \{\delta^{(K+1)}\}}$ be the solution set to the scenario-based \gls{VI} in \eqref{eq:scenario_based_VI} after observing the $(K+1)$-th realization of $\delta$, i.e., the feasible set of the \gls{VI} shrinks to $\mc{X}_{\delta_{K+1}} \coloneqq \mc{X}_{\delta_K} \cap \mc{X}_{\delta^{(K+1)}}$, for some $\delta^{(K+1)} \in \Delta$. We have the following preliminary result.
\begin{lemma}\label{lemma:if_affine}
	For all $K \in \N_0$ and for all $\delta_K \in \Delta^K$, $\Omega_{\delta_{K+1}} = \Omega_{\delta_K} \cap \mc{X}_{\delta^{(K+1)}}$.
\end{lemma}
\begin{pf}
	We split the proof into two different inclusions. Specifically, we first prove (\textrm{i}) that $\Omega_{\delta_K} \cap \mc{X}_{\delta^{(K+1)}} \subseteq \Omega_{\delta_{K+1}}$ , and then (\textrm{ii}) that $\Omega_{\delta_K} \cap \mc{X}_{\delta^{(K+1)}} \supseteq \Omega_{\delta_{K+1}}$.

	\textrm{(i)} We will show that if $x^\star \in \Omega_{\delta_K}$ and $x^\star \in  \mc{X}_{\delta^{(K+1)}}$, then $x^\star \in \Omega_{\delta_{K+1}}$.
	Note that, in view of Standing Assumptions~\ref{standing:sets}--\ref{standing:mapping}, given an arbitrary $K \in \N_0$ and related $K$-multisample $\delta_K \in \Delta^K$, $\mc{X}_{\delta_K}$ is a compact and convex set, as it is finite intersection of convex sets. Then, a vector $x^\star \in \mc{X}_{\delta_K}$ is a solution to VI$(\mc{X}_{\delta_K}, F)$ if and only if $x^\star \in \textrm{argmin}_{y \in \mc{X}_{\delta_K}} \, y^\top F(x^\star)$. Since the uncertain parameter does not affect the mapping $F(\cdot)$, but enters in the constraints only, every sample $\delta^{(K+1)} \in \Delta$ introduces an additional set of convex constraints, i.e., $\mc{X}_{\delta_{K+1}} = \mc{X}_{\delta_K} \cap \mc{X}_{\delta^{(K+1)}} \subseteq \mc{X}_{\delta_K}$, which is compact and convex as well. Thus, it follows immediately that, if $x^\star \in \mc{X}_{\delta^{(K+1)}}$, then $x^\star \in \mc{X}_{\delta_{K+1}}$. Therefore, $x^\star \in \textrm{argmin}_{y \in \mc{X}_{\delta_K} \cap \mc{X}_{\delta^{(K+1)}}} \, y^\top F(x^\ast)$, which by definition implies that $x^\star \in \Omega_{\delta_{K+1}}$.

	\textrm{(ii)} We first prove that, if $x^\star \in \mathrm{relint}(\Omega_{\delta_{K+1}})$, then $x^\star \in \Omega_{\delta_K}$.
	The case where $x^\star \in \textrm{bdry}(\Omega_{\delta_{K+1}})$ will be treated in the sequel.
  Let us recall that, in view of \cite[Cor.~1.6.1]{rockafellar1970convex}, for any given $m$-dimensional convex set $\mc{S}$ in $\R^n$, $m \leq n$, there always exists an affine transformation which carries $\textrm{aff}(\mc{S})$ onto the subspace
	\ifTwoColumn $ \else $$ \fi
	V \coloneqq \{ x = (z_1, \ldots, z_m, z_{m+1}, \ldots, z_n)^\top \in \R^n \mid z_{m+1} = \ldots = z_n = 0 \}.
	\ifTwoColumn $ \else $$ \fi
	Therefore, as closures and relative interiors are preserved under one-to-one affine transformations of $\R^n$ onto itself, we can limit our attention to the case where $\Omega_{\delta_{K+1}}$, and hence $\Omega_{\delta_K}$ (since $\textrm{aff}(\Omega_{\delta_{K+1}}) = \textrm{aff}(\Omega_{\delta_K}) = \textrm{aff}(\Omega_{\delta_0})$ from Standing Assumption~\ref{standing:affinehull}), is $n$-dimensional so that $\textrm{relint}(\Omega_{\delta_{K+1}}) = \textrm{int}(\Omega_{\delta_{K+1}})$.

	 Now, for the sake of contradiction, let $x^\star \in \mc{X}_{\delta_K} \cap \mc{X}_{\delta^{(K+1)}}$ be any point such that $x^\star \in \textrm{int}(\Omega_{\delta_{K+1}})$, but $x^\star \notin \Omega_{\delta_K}$. Since $x^\star \in \mc{X}_{\delta_K}$, $x^\star \notin \Omega_{\delta_K}$ implies that there exists some $\bar{y} \in \mc{X}_{\delta_K}$, with $\bar{y} \neq x^\star$, such that the \gls{VI} is not satisfied, i.e., $(\bar{y} - x^\star)^\top F(x^\star) < 0$. Given the convexity of the sets involved, there must exist some $\lambda \in (0,1)$ that allows one to construct some $\tilde{y} = \lambda x^\star + (1-\lambda) \bar{y}$ such that $\tilde{y} \in \mc{X}_{\delta_K} \cap \mc{X}_{\delta^{(K+1)}}$, but $\tilde{y} \notin \textrm{int}(\Omega_{\delta_{K+1}})$ (see Fig.~\ref{fig:proof_scheme} for a graphical representation). Therefore, since $x^\star \in \textrm{int}(\Omega_{\delta_{K+1}})$, it shall satisfy $(\tilde{y} - x^\star)^\top F(x^\star) \geq 0$, which leads to $(1-\lambda) (\bar{y} - x^\star)^\top F(x^\star) \geq 0$ that clearly generates a contradiction, since $(1 - \lambda) > 0$.

	It remains to show the claim for the case where $x^\star \in \textrm{bdry}(\Omega_{\delta_{K+1}})$. Notice that, since $\textrm{relint}(\Omega_{\delta_{K+1}}) \neq \emptyset$ as $\Omega_{\delta_{K+1}}$ is nonempty, and since the involved sets are closed and convex, for any $x^\star \in \textrm{bdry}(\Omega_{\delta_{K+1}})$ we can always construct a convergent sequence of points $\{x_t\}_{t \in \N}$ such that, for all $t \in \N$, $x_t \in \textrm{relint}(\Omega_{\delta_{K+1}}) \subseteq \Omega_{\delta_K}$, and $\{x_t\}_{t \in \N} \to x^\star$, implying that $x^\star \in \Omega_{\delta_K}$.
	Specifically, given any $\bar{x} \in \textrm{relint}(\Omega_{\delta_{K+1}})$, in view of \cite[Th.~6.1]{rockafellar1970convex}, for all $t \geq 1$, any term of the sequence $x_t \coloneqq \tfrac{1}{t} \bar{x} +  (1 - \tfrac{1}{t}) x^\star \in \Omega_{\delta_{K}} \cap \mc{X}_{\delta^{(K+1)}}$ belongs to $\textrm{relint}(\Omega_{\delta_{K+1}})$.
	Therefore, we have the inclusion $\Omega_{\delta_{K+1}} \subseteq \Omega_{\delta_K}$ as desired.
\end{pf}
\begin{figure}[t]
	\centering
	\ifTwoColumn
		\includegraphics[trim = 0 1cm 0 1cm,
										 width=1.0\columnwidth]{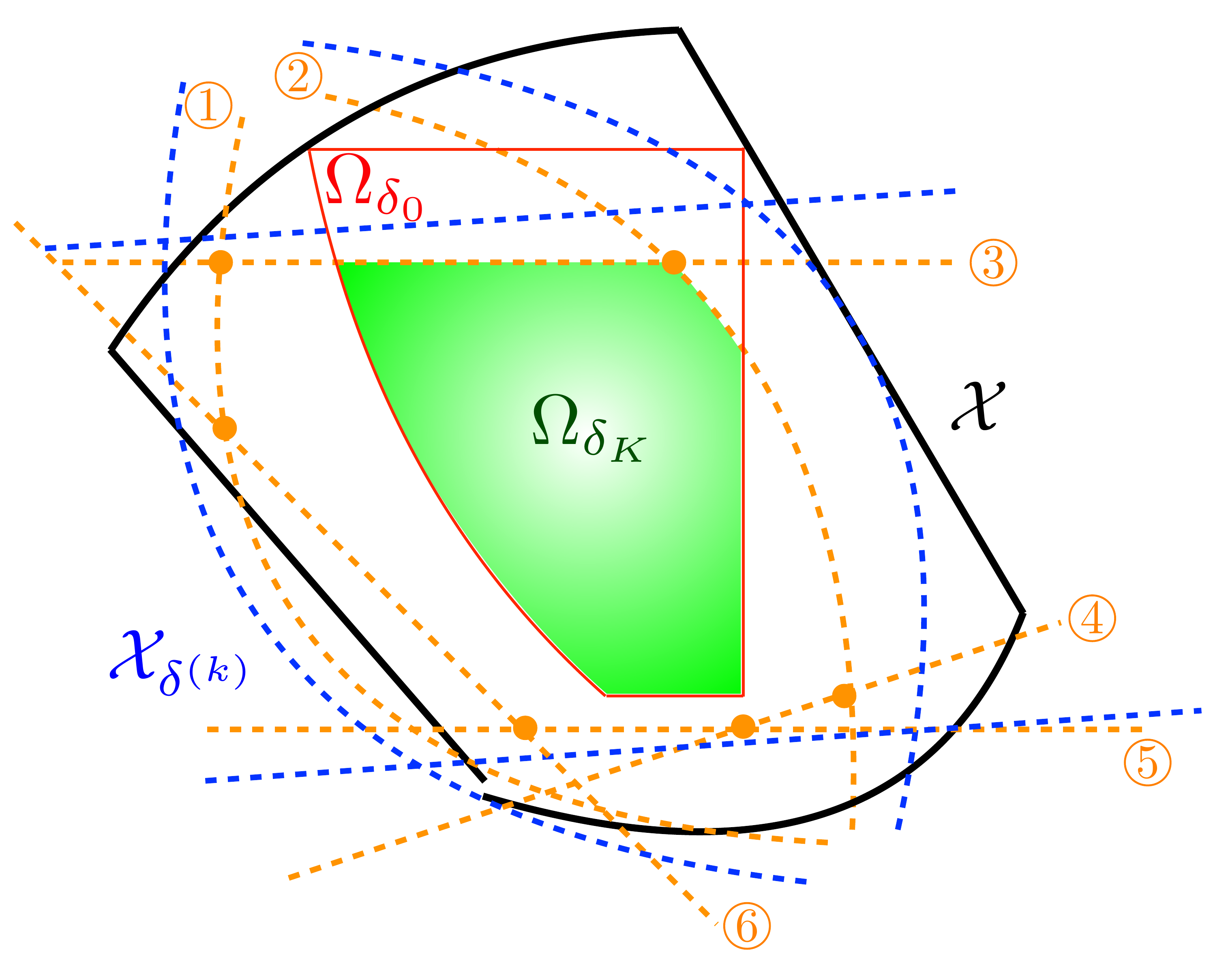}
	\else
		\includegraphics[trim = 0 1cm 0 1cm,
										 width=0.5\columnwidth]{scenario_GNEP}
	\fi
	\caption{The solution set to VI$(\mc{X}_{\delta_K}, F)$, $\Omega_{\delta_K}$ (green region), can be ``shaped'' by the set of constraints, $\mc{X}_{\delta^{(i)}}$, $i \in \mc{K}$ (dashed blue lines). According to Definition~\ref{def:support_sub}, the number of support subsamples for $\delta_K$ \gls{wrt} to $\Omega_{\delta_K}$ is, in general, smaller compared to $|\mc{X}_{\delta_K}|$ (dashed orange lines, whose intersections are defined by orange dots).}
	\label{fig:scenario_based_VI}
\end{figure}
Note that a direct consequence of Lemma~\ref{lemma:if_affine} is that $\Theta_0 \eqqcolon \Omega_{\delta_0} \supseteq \Omega_{\delta_1} \supseteq \ldots \supseteq \Omega_{\delta_K} \eqqcolon \Theta_K(\delta_K)$.
Moreover, the intrinsic consistency of the set $\Omega_{\delta_{K}}$ implies that by introducing additional constraints, the effect of the uncertain parameter is to (possibly) shrink the feasible set on which the scenario-based \gls{VI} in \eqref{eq:scenario_based_VI} is defined, and therefore the set of solutions can only shrink, accordingly (see Fig.~\ref{fig:scenario_based_VI} for a graphical interpretation).

\subsection{Proof of Theorem~\ref{th:VI} and discussion}
We are now ready to prove Theorem~\ref{th:VI}.
\begin{pf}
	For any $K \in \N$, $\delta_{K} \in \Delta^K$, we know that $\Omega_{\delta_K}$ is consistent \gls{wrt} the collected scenarios, $\delta_K$. In view of the definition in \eqref{eq:theta_mapping}, indeed, we have that $\Omega_{\delta_K} \subseteq \cap_{i \in \mc{K}} \mc{X}_{\delta^{(i)}} \cap \mc{X}$, which implies that $\Omega_{\delta_K} \subseteq \mc{X}_{\delta^{(i)}}$, for all $i \in \mc{K}$, thus directly satisfying Definition~\ref{def:cons_set}. Then, we can apply \blue{Lemma~\ref{lemma:campigaratti}} with $\theta^\star_{\delta_K}$ being replaced by $\Omega_{\delta_K}$. We thus have that $\mathbb{P}^K \{\delta_K \in \Delta^K \mid \mathbb{P} \{\delta \in \Delta \mid \Omega_{\delta_K} \not\subseteq \mc{X}_{\delta} \} > \varepsilon(s_K) \} \leq \beta$. However, by Lemma~\ref{lemma:if_affine}, $\Omega_\delta = \Omega_{\delta_K} \cap \mc{X}_\delta$. Therefore, $\Omega_{\delta_K} \not\subseteq \mc{X}_\delta$ is equivalent to $\Omega_\delta \neq \Omega_{\delta_K}$, and since the set of solutions can only shrink once a new scenario is added, this is in turn equivalent to $\Omega_\delta \not\subseteq \Omega_{\delta_K}$. Thus, in view of the definition of $V$ in \eqref{eq:violation_set}, we finally have that $\mathbb{P}^K \{\delta_K \in \Delta^K \mid V(\Omega_{\delta_K}) > \varepsilon(s_K) \} \leq \beta$.
\end{pf}
A more direct expression of the critical parameter  $\varepsilon(\cdot)$ can be obtained by splitting the confidence parameter $\beta$ evenly among the $K$ terms within the summation, i.e.,
\begin{equation}\label{eq:epsilon_analy}
	\varepsilon(h) = \left\{
	\begin{aligned}
	& 1 && \text{ if } \; h = K, \\
	& 1 - \sqrt[K - h]{\frac{\beta}{K \left(\begin{smallmatrix}K\\h\end{smallmatrix}\right)}} && \text{otherwise}.
	\end{aligned}
	\right.
\end{equation}


\blue{
\begin{remark}
	In the case of a non-degenerate \gls{VI}, the bound in \eqref{eq:epsilon} could be improved by means of the \emph{wait-and-judge} analysis in \cite{campi2018wait}. Specifically, in view of \cite[Th.~2]{campi2018wait}, we can replace the expression for $\varepsilon(\cdot)$ in \eqref{eq:epsilon_analy} with $\varepsilon(h) = 1-t(h)$, where $t(h)$ is shown to be the unique solution in $(0,1)$ to
	$$
		\frac{\beta}{K+1} \sum_{m=h}^K {m \choose h} t^{m-h} - {K \choose h} t^{K-h} = 0.
	$$
However, note that the non-degeneracy condition is in general difficult to verify even in convex optimization settings, a challenge that becomes more involved for \glspl{VI}.
\end{remark}
}

Similarly to $\Upsilon(\cdot)$, let us suppose to have available an algorithm that allows us to compute a support subsamples for $\delta_K$ associated with the feasible set $\mc{X}_{\delta_K}$, i.e., the subset of samples such that by using only this subset leads to $\mc{X}_{\delta_K}$ (as opposed to $\Omega_{\delta_K}$, if $\Upsilon(\cdot)$ is employed).
By comparing the bound in \eqref{eq:prob_feas_boud} with the certificates in \cite{pantazis2020aposteriori}, we provide an upper bound for $V(\Omega_{\delta_K})$.
\begin{proposition}\label{prop:better_performance}
	Given any $K \in \N_0$ and $\delta_K \in \Delta^K$, let $s_K$ and $v_K$ be the cardinality of the support subsample for $\delta_K$, evaluated \gls{wrt} $\Omega_{\delta_K}$ and $\mc{X}_{\delta_K}$, respectively. Then, $s_K \leq v_K$.
\end{proposition}
\begin{pf}
	For every $K \in \N_0$ and $\delta_K \in \Delta^K$, Definition~\ref{def:support_sub} suggests that some sample $\delta^{(k)}$ is of support for $\delta_K$ \gls{wrt} $\mc{X}_{\delta_K}$ if $\mc{X}_{\delta^{(k)}}$ is active on $\textrm{bdry}(\mc{X}_{\delta_K})$, i.e., $\textrm{bdry}(\mc{X}_{\delta^{(k)}}) \cap \textrm{bdry}(\mc{X}_{\delta_K}) \neq \emptyset$. On the other hand, $\delta^{(k)}$ is of support \gls{wrt}  $\Omega_{\delta_K}$ if  $\textrm{bdry}(\mc{X}_{\delta^{(k)}}) \cap \Omega_{\delta_K} \neq \emptyset$ (see Fig.~\ref{fig:scenario_based_VI} for a graphic illustration). Since $\Omega_{\delta_K} \subseteq \mc{X}_{\delta_K} \coloneqq \cap_{k \in \mc{K}} \mc{X}_{\delta^{(k)}} \cap \mc{X}$, those samples that are of support for $\delta_K$ \gls{wrt} $\Omega_{\delta_K}$, are necessarily of support \gls{wrt} $\mc{X}_{\delta_K}$, but not vice-versa, and hence $s_K \leq v_K$.
\end{pf}
Under Proposition~\ref{prop:better_performance}, Theorem~\ref{th:VI} improves over \cite{pantazis2020aposteriori}, where $V(\Omega_{\delta_K}) > \varepsilon(v_K)$ was claimed with confidence at most $\beta$. The latter is since $\varepsilon(s_K) \leq \varepsilon(v_K)$ as $\varepsilon(\cdot)$ is non-decreasing. Moreover, within the set-oriented framework proposed in \S 2, as evident from \eqref{eq:prob_feas_boud}, to bound the feasibility risk \eqref{eq:violation_set} of the entire set of solutions $\Omega_{\delta_K}$, one does not need an explicit characterization of $\Omega_{\delta_K}$, namely some mapping $\Theta_K(\cdot)$, but rather the number of support subsamples $s_K$, computed through an algorithm $\Upsilon(\cdot)$. We investigate the computation of $s_K$ in the next section.

\section{\blue{The case of affine constraints:} computation of the support subsample}
In this section we \blue{first propose an iterative procedure that, in case of affine constraints, allows one to compute a support subsample for the unknown $\Omega_{\delta_K}$. We also discuss the related computational complexity}.
\blue{\subsection{Support subsample computation}}
\blue{We start by noting that} the bound in \eqref{eq:prob_feas_boud} depends only on the support subsample that characterizes the solution set $\Omega_{\delta_K}$, and not on its actual shape.
Except in some particular cases, e.g., monotone affine mapping $F(\cdot)$ (see \cite[Th.~2.4.15]{facchinei2007finite}), an explicit representation of $\Omega_{\delta_K}$ is either unavailable, or difficult to compute in advance. \blue{However, the general setting considered so far, i.e., pseudomonotone mapping $F$ and convex constraint set $\mc{X}_{\delta_{K}}$, for any $\delta_{K} \in \Delta^K$, poses several challenges in designing an efficient procedure to compute the number of support subsamples \gls{wrt} $\Omega_{\delta_K}$.}  We therefore \blue{introduce} the following \blue{additional} assumption that restricts attention to the class of linearly constrained\blue{, pseudomonotone} \glspl{VI}.

\begin{assumption}\label{ass:polyhedral}
	Let $\mc{X} \coloneqq \{\bs{x} \in \R^n \mid C \bs{x} \leq d\}$, $C \in \R^{m \times n}$ and $d \in \R^m$, with $\rank(C) = n$, and, for all $\delta \in \Delta$, $\mc{X}_\delta \coloneqq \{\bs{x} \in \R^n \mid A(\delta) \bs{x} \leq b(\delta) \}$, $A : \Delta \to \R^{p \times n}$ and $b : \Delta \to \R^{p}$.
\end{assumption}

\begin{algorithm}[!t]
	\caption{Computation of the cardinality of the support subsample \gls{wrt} $\Omega_{\delta_K}$}\label{alg:support_poly}
	\smallskip
	\textbf{Initialization:}
	\smallskip
	\begin{itemize}\setlength{\itemindent}{-.3cm}
		\item[] Set $s_K \coloneqq 0$, identify
		%
		$\mc{A}_K \coloneqq \{k \in \mc{K} \mid \textrm{bdry}(\mc{X}_{\delta^{(k)}}) \cap \textrm{bdry}(\mc{X}_{\delta_K}) \neq \emptyset \}$
		%
		%
	\smallskip
	\end{itemize}
	\textbf{Iteration $(i \in \mc{A}_K)$:}
	\begin{itemize}\setlength{\itemindent}{.5cm}
		\smallskip
		\item[(\texttt{S1})] Run $\Phi(\delta_{K, i})$ to solve VI$(\mc{X}_{\delta_K} \cap \textrm{bdry}(\mc{X}_{\delta^{(i)}}),F)$
		\smallskip\smallskip
		\item[(\texttt{S2})] If $\Phi(\delta_{K, i}) \neq \emptyset$, set $s_K \coloneqq s_K + 1$
	\end{itemize}
\end{algorithm}

Then, given any $K$-multisample $\delta_K$, let $\Phi : \Delta^K \rightrightarrows\Omega_{\delta_K}$ be any mapping that returns a (set of) solution(s) to VI$(\mc{X}_{\delta_K}, F)$.
The procedure $\Phi(\cdot)$ is run in (\texttt{S1}) to verify whether (at least) one solution to the \gls{VI} with constraints involving $\mc{X}_{\delta_K} \cap \textrm{bdry}(\mc{X}_{\delta^{(i)}})$ exists, thus increasing the counter $s_K$ in case of affirmative answer in (\texttt{S2}). The next result states that, even without knowing the set $\Omega_{\delta_K}$, Algorithm~\ref{alg:support_poly} returns the cardinality of a support subsample for $\delta_K$ \gls{wrt} $\Omega_{\delta_K}$.

\begin{proposition}\label{prop:irreducible_set}
	Let Assumption~\ref{ass:polyhedral} hold true. Given any $K \in \N$ and $\delta_K \in \Delta^K$, Algorithm~\ref{alg:support_poly} returns $s^\star_K$, the cardinality of a support subsample $\delta_K$ \gls{wrt} $\Omega_{\delta_K}$.
\end{proposition}
\begin{pf}
	First note that, in view of Assumption~\ref{ass:polyhedral}, $\mc{A}_K$ forms a support subsample for $\delta_K$ \gls{wrt} $\mc{X}_{\delta_K}$.  Then, by following the considerations adopted within the proof of Proposition~\ref{prop:better_performance}, every $\delta^{(k)}$, $k \in \mc{A}_K$, is of support also  \gls{wrt} to $\Omega_{\delta_K}$ if and only if $\textrm{bdry}(\mc{X}_{\delta^{(k)}}) \cap \Omega_{\delta_K} \neq \emptyset$. To determine this, it is sufficient to compute a solution (if any) on the active region of $\mc{X}_{\delta_K}$ associated with $\mc{X}_{\delta^{(k)}}$. Then, $s_K$ increases only if $\Phi(\delta_{K, k}) \neq \emptyset$, excluding all those samples such that $\mc{X}_{\delta^{(k)}}$ does not intersect $\Omega_{\delta_K}$.
\end{pf}

\begin{remark}
	Algorithm~\ref{alg:support_poly} requires one to run $|\mc{A}_K|$-times the adopted solution algorithm $\Phi(\delta_K)$, with $|\mc{A}_K| \leq K$. This clearly improves \gls{wrt} the greedy algorithms proposed in $\mathrm{\text{\cite[\S II]{campi2018general}}}$ and $\mathrm{\text{\cite[\S III]{paccagnan2019scenario}}}$, which would require one to run $\Phi(\delta_K)$ $K$-times. \blue{On the other hand, we remark that the greedy algorithm applies more generally
	and, in particular, if one assumes uniqueness of the solution as in \cite{campi2018general,paccagnan2019scenario}, it can be employed (not necessarily only in the case of affine constraints) as long as this solution is computable.}
\end{remark}
 \blue{We remark that} the design of an analogous procedure to Algorithm~\ref{alg:support_poly} involving general convex constraints constitutes a topic of current investigation\blue{, as discussed below}.
\blue{\subsection{Computational aspects}}
\blue{From a computational point of view,} we note that Assumption~\ref{ass:polyhedral} is needed for two main reasons:
\begin{enumerate}[leftmargin=*,label=\roman*)]
   \item Evaluating a solution to the \gls{VI} on the boundary of a convex set, i.e., (\texttt{S1}), may require solution of a \gls{VI} defined over a nonconvex domain. Unlike the convex case, the literature on solution algorithms with suitable convergence guarantees or performance for the nonconvex case is not extensive.
	 \item The initialization step becomes trivial, since it requires one to identify the minimal number of active constraints, a task that is closely related to enumerating the number of facets of the polytope $\mc{X}_{\delta_K}$.
 \end{enumerate}

\blue{While item i) prevents us from trivially extending Algorithm~\ref{alg:support_poly} to the case of general convex constraints, item ii) can be equivalently seen as a problem of removing redundant halfspaces. Specifically, in view of the affine structure of both $\mc{X}_\delta$ and $\mc{X}$, the ``offline'' initialization step amounts to solve a family of \glspl{LP}, a task that can be efficiently accomplished in polynomial time by means of available solvers. In fact, given any $K$-multisample $\delta_{K} \in \Delta^K$, the convex polytope $\mc{X}_{\delta_{K}}$ is described by the system of $(m + Kp)$ linear inequalities $A_{\delta_{K}} x \leq b_{\delta_{K}}$, with $A_{\delta_{K}} \coloneqq \col(C, \col((A(\delta_k))_{k \in \mc{K}}))$ and $b_{\delta_{K}} \coloneqq \col(d, \col((b(\delta_k))_{k \in \mc{K}}))$. Here, let $a_i^\top$ (respectively, $b_i$) be the $i$-th row ($i$-th element) of $A_{\delta_{K}}$ ($b_{\delta_{K}}$), and let $\mc{L} \coloneqq \{1, \ldots, m + Kp\}$ be the associated set of row indices. Then, it turns out that a particular $i \in \mc{L}$ is not redundant for $A_{\delta_{K}} x \leq b_{\delta_{K}}$ if and only if the optimal value of the following \gls{LP}
$$
\left\{
\begin{aligned}
&\underset{x \in \R^n}{\textrm{max}} & & a_i^\top x\\
&\hspace{0cm}\textrm{ s.t. } & & a_j^\top x \leq b_j, \, \textrm{ for all } j \in \mc{L}\setminus\{i\},\\
&&& a_i^\top x \leq b_i + 1,
\end{aligned}
\right.
$$
is strictly greater than $b_i$ (notice that the $i$-th constraint has been relaxed). However, particularly when the dimension of the \gls{VI} problem $n$ is large, arbitrarily removing a single inequality constraint at time might prove highly inefficient, or even prohibitive for large amounts of data. For this reason, one may rely on tailored algorithms available in the literature, e.g., \cite{avis1992pivoting,bremner1998primal,christof2001decomposition,ziegler2012lectures}. Successively, once the set of constraints that determines $\textrm{bdry}(\mc{X}_{\delta_{K}})$ has been identified, Algorithm~\ref{alg:support_poly} requires one to run some $\Phi(\cdot)$ to verify whether a solution to the \gls{VI} on each facet of the convex polytope $\mc{X}_{\delta_{K}}$ exists (\texttt{S1}). Since the feasible set of VI$(\mc{X}_{\delta_K} \cap \textrm{bdry}(\mc{X}_{\delta^{(i)}}),F)$ is affine, we remark that the computational complexity of solving (\texttt{S1}) $|\mc{A}_K|$-times is directly driven by the class of (pseudo)monotone mapping $F$ we are tackling, as well as by the family of solution algorithms $\Phi(\cdot)$ we are adopting. As mentioned in \S 2.1, indeed, the literature on efficient algorithms to solve VI$(\mc{X}_{\delta_K} \cap \textrm{bdry}(\mc{X}_{\delta^{(i)}}),F)$ is vast, and therefore an a-priori estimate of the computational burden appears hard to quantify. In \S 5.2, we detail the computational effort when applying this procedure on our case study.
}

\section{Application to robust game theory}

In this section we first discuss how robust game theory and \glspl{GNEP} can be modelled as uncertain \glspl{VI} of the form \eqref{eq:unc_VI}, and then we present a numerical case study modelling the charging coordination of a fleet of \glspl{PEV}. Note that an explicit characterization of the set of solutions to the scenario-based \gls{VI}, $\Omega_{\delta_K}$, is unlikely to exist in this case. Therefore, under appropriate assumptions, i.e., a suitable counterpart of Assumption~\ref{ass:polyhedral}, the systematic procedure in Algorithm~\ref{alg:support_poly} allows us to compute the cardinality of the support subsample.

\subsection{Uncertain \glspl{GNEP} and scenario-based formulation}
We consider a finite population of $N$ agents taking part in a noncooperative game. The agents, indexed by the set $\mc{I} \coloneqq \{1, \ldots,N\}$, control a decision vector $x_i \in \R^{n_i}$, which is locally constrained to a deterministic set $\mc{X}_i \subseteq \R^{n_i}$, which may include, e.g., operational and dynamic constraints over a certain prediction/control horizon.  Each agent aims to minimize a cost function $J_i : \R^n \to \R$, $n \coloneqq \sum_{i \in \mc{I}} n_i$, while satisfying a set of coupling constraints affected by uncertainty $\delta$. The uncertain \gls{GNEP} is hence formalized by the following collection of coupled optimization problems:
\begin{equation}\label{eq:single_prob_orig}
\forall i \in \mc{I} : \left\{
\begin{aligned}
&\underset{x_i \in \mc{X}_i}{\textrm{min}} & & J_i (x_i,  \bs{x}_{-i})\\
&\hspace{.1cm}\textrm{ s.t. } & & (x_i, \bs{x}_{-i}) \in \mc{X}_{\delta}, \, \delta \in \Delta.
\end{aligned}
\right.
\end{equation}
We then introduce the following standard assumptions.
\begin{assumption}\label{ass:game_standing}
	For all $i \in \mc{I}$, for all $\bs{x}_{-i} \in \R^{n - n_i}$, $J_i(\cdot, \bs{x}_{-i})$ is a convex function of class $\mc{C}^1$.
\end{assumption}
\begin{assumption}\label{ass:game_unc}
	For all $i \in \mc{I}$, $\mc{X}_i$ is a nonempty compact, convex set. For all $\delta \in \Delta$, $\mc{X}_\delta$ is a nonempty closed, convex set.
\end{assumption}
Let $\delta_K$ be the $K$-multisample introduced in \S 2, for some $K \in \N_0$. The scenario-based \gls{GNEP} $\Gamma$ is defined as the tuple $\Gamma \coloneqq (\mc{I}, (\mc{X}_i)_{i \in \mc{I}}, (J_i)_{i \in \mc{I}}, \delta_K)$, formally represented by the following family of optimization problems:
\begin{equation}\label{eq:single_prob}
\forall i \in \mc{I} : \left\{
\begin{aligned}
&\underset{x_i \in \mc{X}_i}{\textrm{min}} & & J_i (x_i,  \bs{x}_{-i})\\
&\hspace{.1cm}\textrm{ s.t. } & & (x_i, \bs{x}_{-i}) \in \mc{X}_{\delta^{(k)}}, \, \forall k \in \mc{K}.
\end{aligned}
\right.
\end{equation}
Given its deterministic nature, in view of Assumption~\ref{ass:game_standing} and \ref{ass:game_unc}, the game $\Gamma$ is a jointly convex \glspl{GNEP} \cite[Def.~2]{facchinei2007generalized}. Then, let us define the sets $\mc{X} \coloneqq \prod_{i \in \mc{I}} \mc{X}_i$, $\mc{X}^{\delta_K}_{i}(\bs{x}_{-i}) \coloneqq \{x_i \in \mc{X}_i \mid  (x_i, \bs{x}_{-i}) \in \cap_{k \in \mc{K}} \mc{X}_{\delta^{(k)}} \}$, and $\mc{X}_{\delta_K} \coloneqq \{\bs{x} \in \mc{X} \mid \bs{x} \in \cap_{k \in \mc{K}} \mc{X}_{\delta^{(k)}} \} \subseteq \mc{X}$. We recall now the following key notion of a Nash equilibrium for $\Gamma$:
\begin{definition}\textup{(Generalized Nash Equilibrium)}\label{def:GNE}
	Let $\delta_K \in \Delta^K$ be a given $K$-multisample. The collective strategy $\bs{x}^\star \in \mc{X}_{\delta_K}$ is a $\mathrm{\gls{GNE}}$ of the scenario-based $\mathrm{\gls{GNEP}}$ $\Gamma$ in \eqref{eq:single_prob} if, for all $i \in \mc{I}$,
	$$
		J_i (x^\star_i,  \bs{x}^\star_{-i}) \leq \underset{y_i \in \mc{X}^{\delta_K}_{i}(\bs{x}^\star_{-i})}{\mathrm{min}} \, J_i (y_i,  \bs{x}^\star_{-i}).
	$$
\end{definition}
Clearly, given the dependency on the set of $K$ realizations $\delta_K$, any equilibrium of $\Gamma$ is a random variable itself.

According to Definition~\ref{def:GNE}, a popular subset of \gls{GNE} of a \gls{GNEP} is the one of \glspl{v-GNE} \cite{facchinei2007generalized}, which coincides with the set of collective strategies that solve the \gls{VI} associated with the \gls{GNEP} in \eqref{eq:single_prob_orig}.
 \blue{Specifically, this type of equilibrium problem has certain advantageous structural properties and can be modelled as an uncertain \gls{VI} problem of the type we have considered in this paper. Moreover, it is significant per se as it provides ``larger social stability'' and ``economic fairness'' \cite[\S 5]{cavazzuti2002nash},\cite[Th.~4.8]{facchinei2007generalized}. Thus}, by defining the game mapping $F:\R^n \to \R^n$ as $F(\bs{x}) \coloneqq \col((\nabla_{x_i} J_i(x_i, \bs{x}_{-i}))_{i \in \mc{I}})$, \blue{we formally introduce the class of \glspl{v-GNE} as follows.
\begin{definition}\textup{(Variational Generalized Nash Equilibrium) \cite[Def.~3]{facchinei2007generalized}}
	Let $\delta_K \in \Delta^K$ be a given $K$-multisample defining a jointly convex \gls{GNEP}. A \emph{\gls{v-GNE}} is any solution $\bs{x}^\star \in \mc{X}_{\delta_K}$ to the \gls{GNEP} $\Gamma$ that is also a solution to VI$(\mc{X}_{\delta_{K}}, F)$.
\end{definition}
In summary,} any vector that solves VI$(\mc{X} \cap \mc{X}_\delta, F)$, $\delta \in \Delta$, belongs to the set of \gls{v-GNE} of the \gls{GNEP} in \eqref{eq:single_prob_orig}. Thus, analogously to \eqref{eq:sol_set}, the set of \gls{v-GNE} of the scenario-based \gls{GNEP} in \eqref{eq:single_prob} is
\begin{equation}\label{eq:vGNE_set}
\Omega_{\delta_K} \! \coloneqq \! \{\bs{x} \in \mc{X}_{\delta_K} \! \mid \! (\bs{y} - \bs{x})^\top F(\bs{x}) \geq 0, \; \forall \bs{y} \in \mc{X}_{\delta_K} \}.
\end{equation}
Assume that the set of \gls{GNE} coincides with the one of \gls{v-GNE}. To invoke Lemma~\ref{lemma:solution_compact}, we assume the following:
\begin{assumption}\label{ass:monotone_game}
	The game mapping $F(\cdot)$ is pseudomonotone.
\end{assumption}
Thus, in the spirit of Definition~\ref{def:violation_set}, by relying on the $K$-multisample and the associated scenario-based \gls{GNEP}, we can now employ Theorem~\ref{th:VI} to provide a-posteriori feasibility certificates to the equilibrium set $\Omega_{\delta_K}$ in \eqref{eq:vGNE_set}.

\subsection{Case study: Charging coordination of plug-in electric vehicles}
The problem of coordinating the day-ahead charging of a fleet of \glspl{PEV}, originally introduced in \cite{ma2011decentralized}, can be modelled as a noncooperative \gls{GNEP} \cite{deori2018price,9030152}. Specifically, in the spirit of the previous section, for each \gls{PEV} $i \in \mc{I}$, we consider a discrete-time linear dynamical system $s_i(t+1) = s_i(t) + b_i x_i(t)$, $t \in \N$, where $s_i \in [0,1]$ is the \gls{SoC}, i.e., $s_i = 1$ represents a fully charged battery, while $s_i = 0$ a completely discharged one; $x_i(t) \in [0,1]$ is the charging control input at the specific time instant $t$, and $b_i > 0$ denotes the charging efficiency.
According to a desired level of \gls{SoC} that has to be achieved, the goal of each \gls{PEV} is to acquire (at least) a charge amount $\gamma_i$ within a finite charging horizon $T \in \N$, in order to satisfy the charging constraint $\sum_{t \in \mc{T}} x_i(t) = \bsone_T^\top x_i \geq \gamma_i$, with $\mc{T} \coloneqq \{0, \ldots, T-1\}$ and $x_i \coloneqq \col((x_i(t))_{t \in \mc{T}}) \in \R^T$, while, minimizing its charging cost, $J_i(x_i, \bs{x}_{-i}) \coloneqq p(\bs{x})^\top x_i$. Here, $p : \R_{\geq 0}^T \to \R^T_{\geq 0}$, denotes the electricity price function over the charging horizon, which for simplicity we assume to be affine in the aggregate demand of energy associated with the set of \glspl{PEV}, namely $p(\bs{x}) \coloneqq \alpha \sigma(\bs{x}) + \eta$, with $\sigma(\bs{x}) \coloneqq \sum_{j \in \mc{I}} x_j \in \R^T$, for some $\alpha > 0$ and $\eta \in \R^{T}_{\geq 0}$. Moreover, due to the intrinsic limitations of the grid capacity $d_{\textrm{max}} > 0$, we assume that the amount of energy required in each single time period by both the \glspl{PEV} and non-\gls{PEV} loads should not be greater than $d_{\textrm{max}}$. This translates into a constraint on the \glspl{PEV} total demand, i.e., $d(t) + \sum_{j \in \mc{I}} x_j(t) \in [0, d_{\textrm{max}}]$, for all $t \in \mc{T}$.

The inflexible non-\gls{PEV} demand $d \in \R_{\geq 0}^T$ is subject to uncertainty and therefore is modelled as $d \coloneqq d_{\textrm{nom}} + \delta$. Here, $d_{\textrm{nom}} \in \R^T_{\geq 0}$ is the nominal non-\gls{PEV} daily energy demand, which can be extracted, e.g., from data (see \cite{natgrid} for typical daily energy profiles in the UK), while $\delta$ is a random variable that follows a uniform probability distribution on $\Delta \subseteq \R^T$. The (uncertain) \gls{GNEP} coincides with the following collection of optimization problems
\begin{equation}\label{eq:PEV_GNEP}
	\forall i \!\in\! \mc{I} \!:\! \left\{
	\begin{aligned}
		&\underset{x_i \in [0,1]^T}{\textrm{min}} & & (\alpha \sigma(\bs{x}) + \eta)^\top x_i\\
		&\hspace{.3cm}\textrm{ s.t. } & & (d_{\textrm{nom}} \!+\! \delta) \!+\! \sigma(\bs{x}) \leq \bsone_{T} d_{\textrm{max}}, \, \forall \delta \!\in\! \Delta,\\
		&&& A_i x_i \leq c_i,\\
	\end{aligned}
	\right.
\end{equation}
where $A_i \coloneqq \col(-B_i, B_i, -\bsone_T^\top) \in \R^{(2T + 1) \times T}$, $B_i \in \R^{T \times T}$ is matrix with all entries in the lower triangular part equal to $b_i$, $c_i \coloneqq \col(\bsone_T s_i(0), \bsone_T (1 - s_i(0)), -\gamma_i) \in \R^{2T + 1}$, and $s_i(0) \in [0,1]$ is a given initial \gls{SoC}. Problem \eqref{eq:PEV_GNEP} is in the form of \eqref{eq:single_prob_orig}.
We note that the game mapping $F(\bs{x}) \coloneqq \col(\nabla_{x_i}((\alpha \sigma(\bs{x}) + \eta)^\top x_i)_{i \in \mc{I}})$, which allows us to define the \gls{VI} whose solution set determines the \glspl{v-GNE} of the game, turns out to be affine in $\bs{x}$. Specifically, $F(\bs{x}) = M \bs{x} + q$, where $M \in \R^{NT \times NT}$ has entries all equal to $\alpha$, while $q \coloneqq  \bsone_N \otimes \eta \in \R^{NT}$.
Note that, for any $\alpha > 0$, $F(\cdot)$ is a monotone mapping (or, equivalently, $M + M^\top \succcurlyeq 0$).

Thus, based on $K$ observations of historical data, the \gls{GNEP} in \eqref{eq:PEV_GNEP} admits a scenario-based counterpart as in \eqref{eq:single_prob}, for which we quantify the robustness of $\Omega_{\delta_K}$, the solution set to VI$(\mc{X}_{\delta_K}, F)$. Here, $\mc{X}_{\delta_K} \coloneqq \prod_{i \in \mc{I}} \mc{X}_i \cap_{k \in \mc{K}} \mc{X}_{\delta^{(k)}}$, with $\mc{X}_i \coloneqq \{x_i \in [0,1]^T \mid A_i x_i \leq c_i \}$, and $\mc{X}_{\delta^{(k)}} \coloneqq \{\bs{x} \in \R^{NT} \mid (d_{\textrm{nom}} + \delta^{(k)}) + \sigma(\bs{x}) \leq \bsone_{T} d_{\textrm{max}}\}$, $k \in \mc{K}$.

\subsubsection{Numerical simulations}
\begin{table}[tb]
	\caption{Simulation parameters}
	\label{tab:sim_param}
	\begin{center}
		\begin{tabular}{ p{1cm} p{3.5cm}p{2.7cm}  }
			\toprule
			Name &  Description & Value \\
			\midrule
			$N$ &  \glspl{PEV} number   &$20$  \\
			$T$ &  Time intervals & $24$ \\
			$b_i$  & Charging efficiency & $[0.075, 0.25]$ \\
			$s_i(0)$ &Initial SoC of battery & $[0.1, 0.4]$\\
			$s_i(T)$ &Desired SoC of battery & $[0.7, 1]$\\
			$\gamma_i$& Required charge amount & $[1.62, 7.49]$ \\
			$\alpha$ & Inverse of price elasticity  & $0.01$ \\
			$\eta$  & Baseline price & $\bs{0}_T$ \\
			$d_{\textrm{nom}}$ &Non-\gls{PEV} demand & Average over $10^2$ daily profiles \cite{natgrid}  \\
			$d_{\textrm{max}}$ & Grid power capacity & $2 \cdot \underset{t \in \mc{T}}{\textrm{max}} \; d_{\textrm{nom}}(t)$ \\
			$\Delta$ & Uncertainty support  & $ d_{\textrm{nom}} \cdot [-0.05, 0.05]$\\
			\bottomrule
		\end{tabular}
	\end{center}
\end{table}


\begin{figure}
	\centering
	\ifTwoColumn
	\includegraphics[width=1.0\columnwidth]{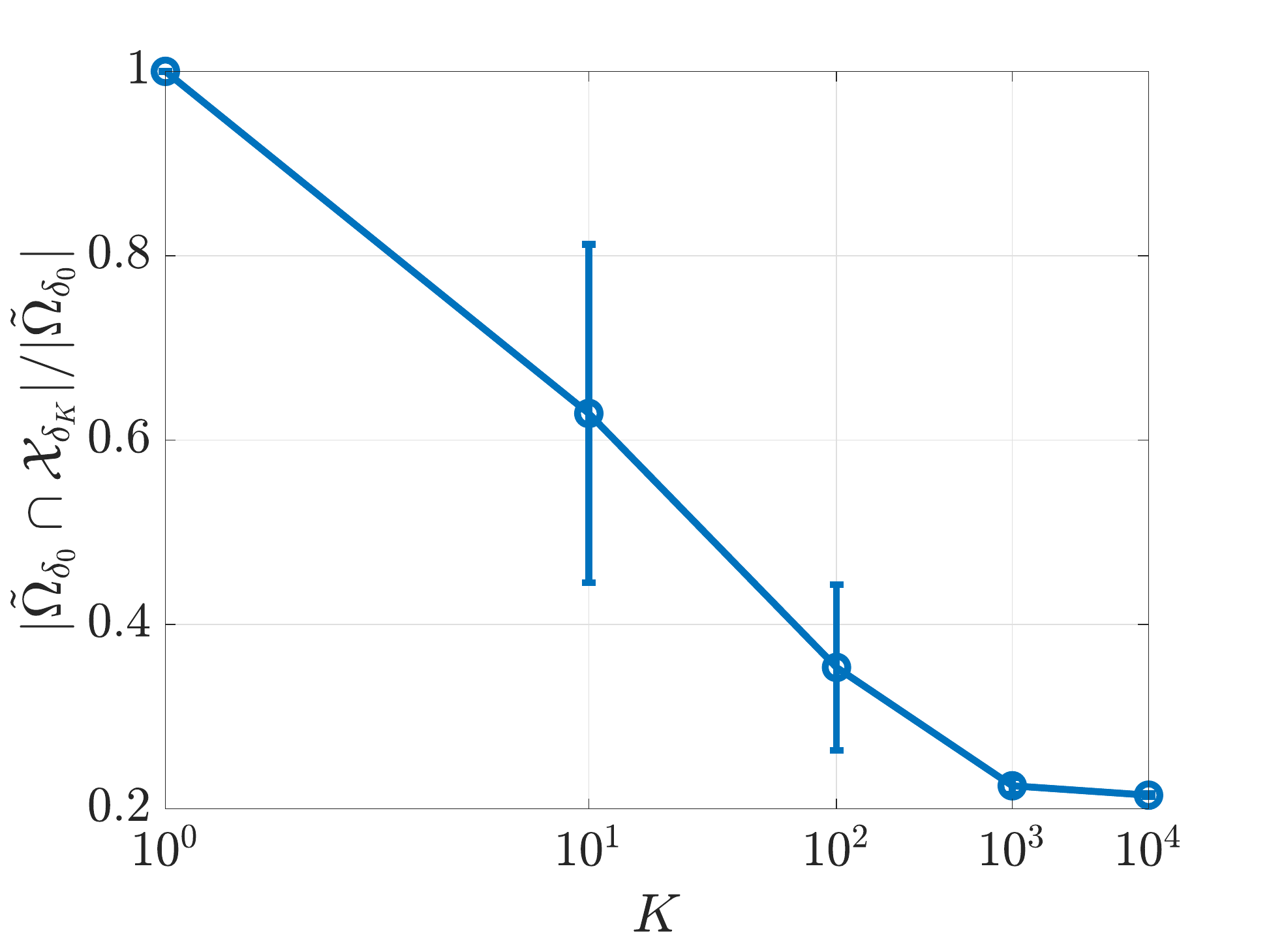}
	\else
	\includegraphics[width=0.5\columnwidth]{reduction-eps-converted-to.pdf}
	\fi
	\caption{Number of solutions contained in $\tilde{\Omega}_{\delta_K} = \tilde{\Omega}_{\delta_0} \cap \mc{X}_{\delta_K}$, normalized with the ones lying in $\tilde{\Omega}_{\delta_0}$, as a function of the number of samples $K$. The solid blue line represents the average of $|\tilde{\Omega}_{\delta_0} \cap \mc{X}_{\delta_K}|/|\tilde{\Omega}_{\delta_0}|$ over $100$ numerical experiments, while the vertical blue lines the standard deviation.}\label{fig:reduction}
\end{figure}

\begin{figure}
	\centering
	\ifTwoColumn
	\includegraphics[width=1.0\columnwidth]{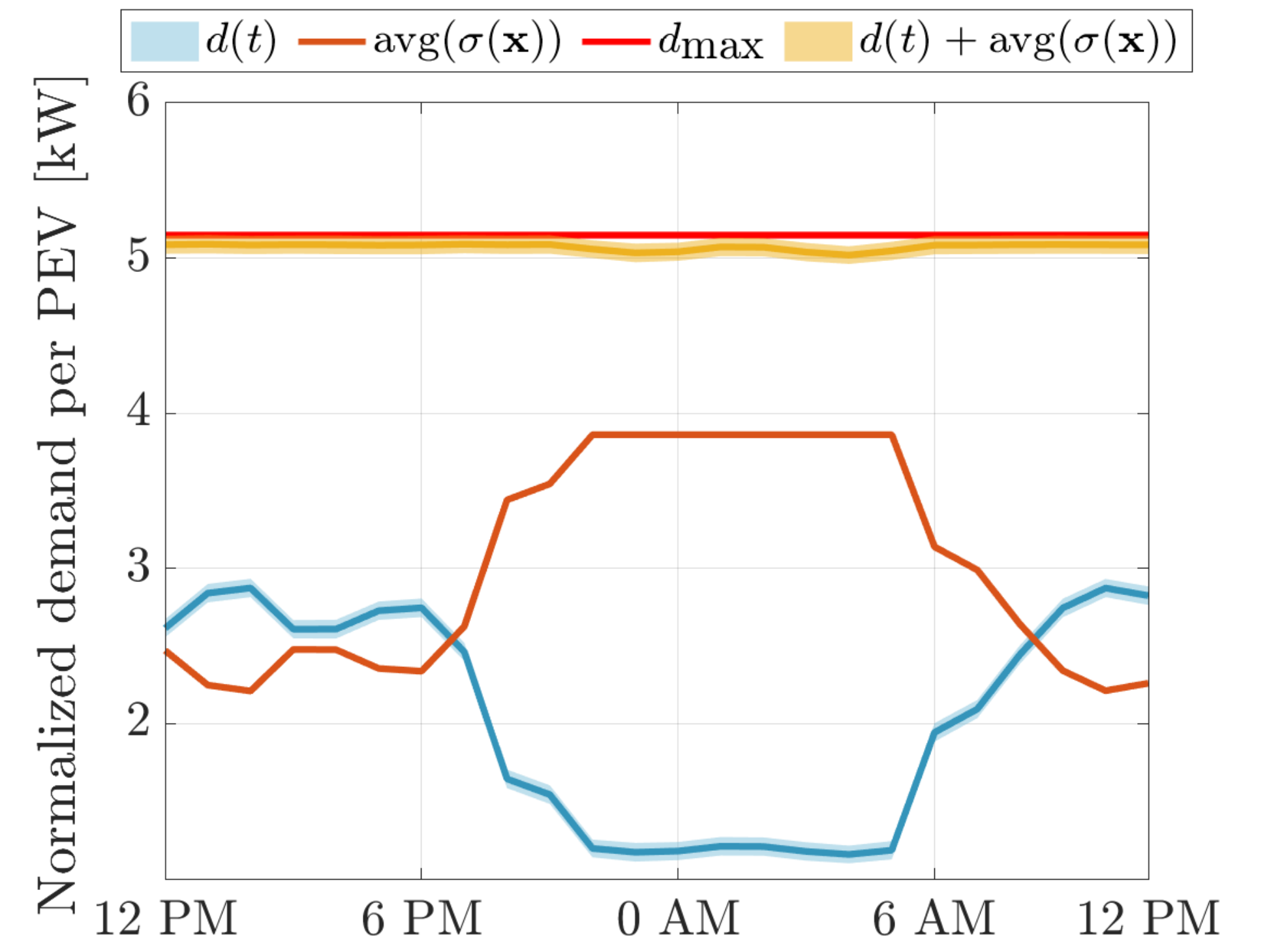}
	\else
	\includegraphics[width=0.5\columnwidth]{avg_agg-eps-converted-to.pdf}
	\fi
	\caption{Average behaviour of the fleet of \glspl{PEV}, computed across the estimated set of solutions $\tilde{\Omega}_{\delta_K}$ after observing $10^3$ realizations of the uncertainty. The overall demand, affected by the uncertainty $\delta$, meets the grid capacity limitations.}\label{fig:avg_agg}
\end{figure}

\blue{To numerically test the theoretical results provided in the paper, in this section we fit a multidimensional Gaussian distribution to $10^2$ daily profiles from \cite{natgrid} to generate large sets of realistic samples. Thus,} we first support the consistency of $\Omega_{\delta_{K}}$ numerically. Specifically, we estimate $\Omega_{\delta_0}$ by computing $10^3$ different solutions to VI$(\mc{X}_{\delta_0}, F)$, hence obtaining $\tilde{\Omega}_{\delta_0}$, with the numerical parameters reported in Table~\ref{tab:sim_param}. Every solution is computed by means of a typical extragradient algorithm \cite{nguyen2018extragradient}, initialized with a different condition and fixed step size, whose convergence is guaranteed as $F$ is monotone and Lipschitz continuous with constant $\alpha NT$, for any step size $(0, \alpha^{-1}/NT)$. We emphasize that we are interested in computing a set of solutions to \gls{VI}$(\mc{X}_{\delta_K}, F)$, and this motivates us to partially neglect the multi-agent nature of the problem addressed by adopting an extragradient method\footnote[1]{Decentralized equilibrium computation is, indeed, outside the scope of the current paper.}.
 \blue{Moreover, to compute a solution to VI$(\mc{X}_{\delta_0}, F)$ with a precision in norm of $10^{-6}$, the extragradient algorithm in \cite{nguyen2018extragradient} takes around $21.47$[s] on average, resulting in $18$--$24$ iterations. Given the linearity of the constraints, this value is representative for solving (\texttt{S1}) in Algorithm~\ref{alg:support_poly}.}
Thus, as illustrated in Fig.~\ref{fig:reduction}, the average number of solutions contained in $\tilde{\Omega}_{\delta_K}$ over $100$ numerical experiments, normalized \gls{wrt} $\tilde{\Omega}_{\delta_0}$, shrinks as $K$ grows. Note that, in view of the structure of $\Delta$, as the number of samples $K$ increases, the standard deviation of $\delta$ narrows around the average. An example of aggregate behaviour of the fleet of \glspl{PEV} is reported in Fig.~\ref{fig:avg_agg}, where $\textrm{avg}(\cdot)$ returns the average among the solutions lying in $\tilde{\Omega}_{\delta_K}$, estimated after observing $10^3$ realizations of the uncertainty. Note that $\sigma(\bs{x})$ exhibits the so-called ``valley filling'' property, which is desirable since the overall demand has no peaks.
\begin{table}[tb]
	\caption{Robustness certificate \eqref{eq:prob_feas_boud} and empirical violation probability}
	\label{tab:feas_bound}
	\begin{center}
		\begin{tabular}{cccccc}
			\toprule
			$K$ &  $|\mc{A}_K|$ & $s^\star_K$  & $\varepsilon(s^\star_K)$ & $V_{\textrm{max}}(\tilde{\Omega}_{\delta_K})$ & $\textrm{avg}(V_{\textrm{max}}(\tilde{\Omega}_{\delta_K}))$\\
			\midrule
			$10^2$ & $28$ & $4$ & $    0.29$ & $0.018$ & $0.012$\\
			$10^3$ & $381$ & $7$ & $0.06$ & $1.2\cdot10^{-3}$ & $0.8\cdot10^{-3}$\\
			$10^4$ & $469$ & $9$ & $0.01$ & $0.9\cdot10^{-3}$ & $0.4\cdot10^{-3}$\\
			\bottomrule
		\end{tabular}
	\end{center}
\end{table}

For any $K \in \N_0$, the feasible set of the scenario-based counterpart of \eqref{eq:PEV_GNEP} satisfies Assumption~\ref{ass:polyhedral}. Thus, in Table~\ref{tab:feas_bound} we compare the output of the procedure summarized in Algorithm~\ref{alg:support_poly} to compute the cardinality $s^\star_K$ of the support subsample \gls{wrt} $\Omega_{\delta_K}$, for different values of $K$. The bound on the violation probability is then computed by the function $\varepsilon(\cdot)$ in \eqref{eq:epsilon_analy} with $\beta = 10^{-6}$. Note that Algorithm~\ref{alg:support_poly} requires us to run $\Phi(\cdot)$ only $|\mc{A}_K|$-times, which represents a noticeable improvement compared to the greedy algorithm proposed in \cite{campi2018general,paccagnan2019scenario}, which would require running $\Phi(\cdot)$ $K$-times. \blue{On the other hand, the ``offline'' initialization step with $K = 10^{4}$, which translates into $241940$ linear inequalities, takes around $11600$[s] to identify the set of constraints defining $\mc{X}_{\delta_{10^{4}}}$. Finally,} the last two columns in Table~\ref{tab:feas_bound} report both the maximum and the average value of the empirical violation probability of $\tilde{\Omega}_{\delta_0} \cap \mc{X}_{\delta_K}$ computed against $10^2$, $10^3$ and $10^4$ new realizations. The empirical probability, as expected, is lower than the theoretical bound in Theorem~\ref{th:VI}.

\section{Concluding remarks}
The scenario approach paradigm applied to uncertain \glspl{VI} provides a numerically tractable framework to compute solutions with quantifiable robustness properties in a distribution-free fashion. In the specific family of uncertain \glspl{VI} considered, which encompasses a broad class of practical applications belonging to different domains, we are able to evaluate the robustness properties of the entire set of solutions, thereby relaxing the requirement of a unique solution as often imposed in the literature. We have shown that this only requires us to enumerate the active coupling constraints that ``shape'' that set.

Future research directions involve synthesizing algorithms to enumerate the number of support subsamples in a general convex setting, as well as investigating extensions of the approach we have developed to quasi-variational inequalities. This would enable us to incorporate the uncertainty within the mapping defining the \gls{VI}, thus extending the results of \cite{fele2019probabilistic,fele2019probably} to the entire set of solutions to \glspl{VI}.




\balance

\bibliographystyle{plain}
\bibliography{scenario_VI}

\end{document}